\documentclass[11pt]{article}
\usepackage[a4paper,margin=28mm]{geometry}
\usepackage{amsmath,amssymb,amsthm,mathtools}
\usepackage{enumitem,microtype,needspace}
\usepackage[colorlinks=true,linkcolor=blue,citecolor=blue,urlcolor=blue]{hyperref}
\hypersetup{pdftitle={Sharp moduli of continuity for the fractional p-Laplacian at the critical order},pdfauthor={Yiyang Xiao}}
\numberwithin{equation}{section}
\newtheorem{theorem}{Theorem}[section]
\newtheorem{proposition}[theorem]{Proposition}
\newtheorem{lemma}[theorem]{Lemma}
\newtheorem{corollary}[theorem]{Corollary}
\theoremstyle{definition}
\newtheorem{definition}[theorem]{Definition}
\theoremstyle{remark}
\newtheorem{remark}[theorem]{Remark}
\newcommand{\R}{\mathbb R}
\newcommand{\LL}{\mathcal L}
\newcommand{\PV}{\operatorname{PV}}
\newcommand{\norm}[1]{\lVert#1\rVert}
\newcommand{\one}{\mathbf 1}
\DeclareMathOperator{\Lip}{Lip}
\DeclareMathOperator{\Tail}{Tail}
\setlist[enumerate]{label=\textup{(\roman*)},leftmargin=*,itemsep=4pt}
\allowdisplaybreaks[2]

\title{Sharp moduli of continuity for the fractional $p$-Laplacian at~the~critical~order}
\author{Yiyang Xiao \thanks{School of Mathematical Sciences, Beijing Normal University, Beijing 100875, China. Email:
202531130042@mail.bnu.edu.cn}}
\date{}

\begin{document}
\maketitle
\begin{abstract}
At the critical order $sp=p-1$, we determine the sharp interior modulus
of continuity for weak solutions of $(-\Delta_p)^s u=f$ with $p>1$
and locally bounded $f$. If $\omega(t)$ bounds the oscillation of $f$
at distances up to $t$, then the interior modulus of continuity of $u$
at scale $r$ is bounded by
\[
 C r\left(1+\int_r^1\omega(t)\,\frac{\mathrm dt}{t}\right)^{1/(p-1)}.
\]
For each nonzero subadditive modulus $\omega$, we construct a bounded
solution that realizes both the prescribed oscillation of $f$
and this rate for $u$, up to constants, at every sufficiently small scale.
For bounded $f$, the sharp interior modulus of continuity of $u$ is
$r(\log(\mathrm e/r))^{1/(p-1)}$.
The above bound gives local Lipschitz continuity for Dini-continuous $f$,
and for every non-Dini modulus our construction gives a solution with
unbounded difference quotients at the origin.
\end{abstract}


\section{Introduction and main results}\label{sec:introduction}

We consider the fractional $p$-Laplace equation $(-\Delta_p)^s u=f$ in
$B_2$, where $n\ge1$, $1<p<\infty$ and $0<s<1$. We use the normalization
\begin{equation}\label{eq:fracplap}
 (-\Delta_p)^s u(x)=\PV\int_{\R^n}
 \frac{|u(x)-u(y)|^{p-2}(u(x)-u(y))}{|x-y|^{n+sp}}\,\mathrm dy.
\end{equation}
For locally bounded $f$, weak solutions are locally H\"older continuous
with every exponent $0<\gamma<\min\{1,sp/(p-1)\}$.
This follows from the estimates of Brasco, Lindgren and Schikorra
\cite[Theorem~1.4]{BLS} for $p\ge2$ and of Garain and Lindgren
\cite[Theorem~1.2]{GL} for $1<p<2$. Away from $sp=p-1$, Biswas and Topp
\cite[Theorem~1.1]{BT} reached the endpoint, proving local H\"older
continuity with exponent $sp/(p-1)$ when $sp<p-1$ and local Lipschitz
continuity when $sp>p-1$.

At the critical order $sp=p-1$, Biswas and Topp \cite[Theorem~1.2]{BT}
proved, for every $p>1$ and locally bounded $f$, a log-Lipschitz
estimate with some logarithmic exponent $\theta>0$. They also
proved local Lipschitz continuity when the modulus of $f$ is
$O((1+|\log t|)^{-(n+3)})$ as $t\downarrow0$. In the linear case
$p=2$, $s=1/2$, the interior Zygmund estimate of Dipierro, Ros-Oton,
Serra and Valdinoci \cite[Theorem~1.1(a)]{DRSV} yields the modulus
$r\log(\mathrm e/r)$ for globally bounded weak solutions with locally
bounded data (see \cite[Remark~2.38 and Proposition~2.107]{BCD}).
Local Lipschitz continuity can fail even in this linear case \cite[Example~1.6]{BLS}.

Very recently, Giovagnoli, Jesus and Silvestre \cite[Theorem~1.1]{GJS}
proved that, at the critical order and for $p\ge2$, solutions of the
homogeneous equation have locally H\"older continuous gradients.
Giovagnoli and Jesus \cite[Theorem~1.2]{GJ} later proved that, in the
same range, solutions with H\"older continuous right-hand sides also
belong to $C^{1,\alpha}_{\mathrm{loc}}$ for some $\alpha>0$.

Our main result shows that when $sp=p-1$, the difference quotients of $u$
satisfy a sharp upper bound in terms of
$\int_{|x-y|}^1 \omega(f;t)\,\frac{\mathrm dt}{t}$, where $\omega(f;t)$ is the
oscillation of $f$ over pairs of points in $B_{3/2}$ at distance
at most $t$. For each nonzero,
nondecreasing, subadditive modulus $\omega$, we construct in
Theorem~\ref{thm:sharp-modulus} a bounded weak solution for which, at
every sufficiently small scale $r$, the oscillation of $f$ is comparable
to $\omega(r)$ and the difference quotient $(u(re_1)-u(0))/r$ matches
the upper bound up to constants. In particular,
for bounded $f$ the optimal interior modulus of continuity of $u$ at
scale $r$ is $r(\log(\mathrm e/r))^{1/(p-1)}$.

To state our main results, we introduce the tail space
\[
 L^{p-1}_{sp}(\R^n)=\left\{u\in L^{p-1}_{\mathrm{loc}}(\R^n):
 \int_{\R^n}\frac{|u(w)|^{p-1}}{(1+|w|)^{n+sp}}\,\mathrm dw<\infty
 \right\}
\]
and the associated nonlocal tail
\begin{equation}\label{eq:tail}
 \Tail(u;x_0,R)=\left(R^{sp}
 \int_{\R^n\setminus B_R(x_0)}
 \frac{|u(w)|^{p-1}}{|w-x_0|^{n+sp}}\,\mathrm dw
 \right)^{1/(p-1)}.
\end{equation}

\Needspace{7\baselineskip}
\begin{theorem}\label{thm:log}
Let $1<p<\infty$ and $s=(p-1)/p$, and let
$u\in W^{s,p}_{\mathrm{loc}}(B_2)\cap L^{p-1}_{sp}(\R^n)$ be a local weak
solution of $(-\Delta_p)^s u=f$ in $B_2$, with
$f\in L^\infty_{\mathrm{loc}}(B_2)$.
Then, for all distinct $x,y\in B_1$,
\begin{equation}\label{eq:logestimate}
 \begin{split}
 |u(x)-u(y)|\le{}&C\left(\norm{u}_{L^\infty(B_{3/2})}+\Tail(u;0,3/2)
                   +\norm{f}_{L^\infty(B_{3/2})}^{1/(p-1)}\right)\\
 &\qquad\times |x-y|\bigl(1+\bigl|\log |x-y|\bigr|\bigr)^{1/(p-1)},
 \end{split}
\end{equation}
where $C$ depends only on $n$ and $p$.
\end{theorem}

The logarithmic factor in \eqref{eq:logestimate} arises when the
oscillation of $f$ at every scale is bounded by
$2\norm{f}_{L^\infty(B_{3/2})}$. Our main estimate retains the
dependence of this oscillation on the scale. For $0<t\le1$, we denote
the essential modulus of continuity of $f$ on $B_{3/2}$ by
\begin{equation}\label{eq:omega}
 \omega(f;t)=\operatorname*{ess\,sup}
 \bigl\{|f(x)-f(y)| : x,y\in B_{3/2},\ |x-y|\le t\bigr\},
\end{equation}
where the essential supremum is taken with respect to product Lebesgue
measure.

\begin{theorem}\label{thm:logmodulus}
Let $1<p<\infty$ and $s=(p-1)/p$, and let
$u\in W^{s,p}_{\mathrm{loc}}(B_2)\cap L^{p-1}_{sp}(\R^n)$ be a local weak
solution of $(-\Delta_p)^s u=f$ in $B_2$, with
$f\in L^\infty_{\mathrm{loc}}(B_2)$. Then, for all $x,y\in B_1$ with
$0<|x-y|\le1$,
\begin{equation}\label{eq:mainmodulus}
 \begin{split}
 |u(x)-u(y)|\le C|x-y|\Bigg[&\norm{u}_{L^\infty(B_{3/2})}+\Tail(u;0,3/2)\\
 &+\left(\norm{f}_{L^\infty(B_{3/2})}
       +\int_{|x-y|}^1\omega(f;t)\frac{\mathrm dt}{t}
   \right)^{1/(p-1)}\Bigg],
 \end{split}
\end{equation}
where $C$ depends only on $n$ and $p$.
\end{theorem}

Theorem~\ref{thm:sharp-modulus} shows that \eqref{eq:mainmodulus} is
sharp for each nondecreasing, subadditive modulus $\omega$ normalized by
$\omega(1)=1$. It gives a bounded local weak solution for which, at every
sufficiently small scale $r$,
\[
 \omega(f;r)\asymp\omega(r),\qquad
 u(re_1)-u(0)\asymp
 r\left(1+\int_r^1\omega(t)\frac{\mathrm dt}{t}\right)^{1/(p-1)},
\]
where $e_1$ is the first coordinate vector. The datum is continuous
when $\omega(0+)=0$.

If the essential modulus of $f$ is Dini, the integral in
\eqref{eq:mainmodulus} remains bounded as $|x-y|\downarrow0$, yielding
the following interior Lipschitz estimate.

\begin{corollary}\label{cor:dini}
Under the hypotheses of Theorem~\ref{thm:logmodulus}, if
$\int_0^1\omega(f;t)\,\mathrm dt/t<\infty$, then $u$ is Lipschitz in
$B_1$ and
\begin{equation}\label{eq:diniestimate}
 \begin{split}
 [u]_{\Lip(B_1)}\le C\Bigg[&\norm{u}_{L^\infty(B_{3/2})}+\Tail(u;0,3/2)\\
 &+\left(\norm{f}_{L^\infty(B_{3/2})}
       +\int_0^1\omega(f;t)\frac{\mathrm dt}{t}
   \right)^{1/(p-1)}\Bigg].
 \end{split}
\end{equation}
Here $C$ depends only on $n$ and $p$.
\end{corollary}
\begin{proof}
For $0<|x-y|\le1$, bound the integral in \eqref{eq:mainmodulus} by
$\int_0^1\omega(f;t)\,\mathrm dt/t$. For $x,y\in B_1$ with $1<|x-y|<2$,
apply this bound from each point to the midpoint $(x+y)/2$.
\end{proof}

Conversely, for every non-Dini modulus $\omega$ as in
Theorem~\ref{thm:sharp-modulus}, there is a bounded weak solution with
$\omega(f;t)\le\omega(t)$ for $0<t\le1$ such that
\[
 \frac{u(re_1)-u(0)}{r}\longrightarrow\infty
 \qquad\text{as }r\downarrow0.
\]
Thus the Dini condition is the exact condition on a prescribed data
modulus that guarantees local Lipschitz continuity.
Corollary~\ref{cor:sharp-dini} gives the corresponding uniform estimate.
In particular, for data satisfying
$\omega(f;t)\le C(1+|\log t|)^{-\mu}$ with $\mu>0$, the sharp threshold
for local Lipschitz regularity is $\mu>1$ (Remark~\ref{rem:sharp-endpoints}).

Throughout Sections~\ref{sec:preliminaries}--\ref{sec:sharpness} we fix
$1<p<\infty$ and $s=(p-1)/p$. We prove the interior estimates for
translation-invariant operators whose kernels $K$ are even, continuous
away from the origin, and satisfy
\[
 \lambda |z|^{-n-sp}\le K(z)\le\Lambda |z|^{-n-sp}
 \qquad(z\ne0,\ 0<\lambda\le\Lambda<\infty).
\]
Under these assumptions we obtain the bounded-data estimate \eqref{eq:logestimate}.
If, in addition, $k(z)=|z|^{n+sp}K(z)$ is Dini continuous on
$\{|z|\ge1/16\}$, we obtain the refined estimate \eqref{eq:mainmodulus}.

To prove these estimates, we use the nonlocal Ishii--Lions method of Biswas
and Topp \cite{BT}. In the doubling-of-variables argument, we compare the
operator at $\bar x$ and $\bar y$ when $u(x)-u(y)-L\phi(|x-y|)-\psi(x)$
attains a positive maximum at $(\bar x,\bar y)$, where $\psi$ is a
localization function. The cone used in the critical
estimate of \cite{BT} narrows as the contact points approach each other.
We use an annular cone of fixed aperture. Writing
$\bar e=(\bar x-\bar y)/|\bar x-\bar y|$, we pair each displacement $z$ with
its reflection
\[
 Rz=z-2(z\cdot\bar e)\bar e.
\]
The displacements $z$ and $Rz$ have the same component perpendicular to
$\bar e$. On the annular cone, the distances between the points in
$(\bar x+z,\bar y+Rz)$ and $(\bar x+Rz,\bar y+z)$ are
$|\bar x-\bar y|+2z\cdot\bar e$ and $|\bar x-\bar y|-2z\cdot\bar e$,
respectively. Applying maximality to these two pairs and adding the
resulting inequalities gives a one-dimensional second difference of $\phi$. A
truncated monotonicity estimate for $J_p(t)=|t|^{p-2}t$ then gives,
after integration over the annular cone, a lower bound proportional to
\[
 -r\frac{\mathrm d}{\mathrm dr}\bigl[(L\phi'(r))^{p-1}\bigr]
 \qquad\text{at }r=|\bar x-\bar y|.
\]
To make this bound dominate the oscillation of $f$, we choose
$(\phi')^{p-1}$ by integrating a majorant of that oscillation against
$\mathrm dt/t$. This gives the integral and the power $1/(p-1)$ in
\eqref{eq:mainmodulus}.

For sharpness we seek a bounded odd function $v$ on the line whose
increments attain the rate in \eqref{eq:mainmodulus} and whose right-hand
side has the prescribed modulus $\omega$, up to constants. Let $\Omega$
be a smooth concave function comparable to $\omega$. We prescribe the
variation of the slope by
\[
 (v'(r))^{p-1}-(v'(1))^{p-1}
 =\int_r^1 q(t)\,\frac{\mathrm dt}{t},\qquad 0<r\le1.
\]
Taking $q=\Omega$ gives the desired growth of $v$. The difficulty arises
in estimating $(-\Delta_p)^s v(r)$. The bound for displacements larger
than $2r$ contains a term of the form $r\int_r^R q(t)\,\mathrm dt/t^2$.
With $q(t)=\Omega(t)=t$ on $[r,R]$, this term is $r\log(R/r)$, whereas
the desired bound is of order $r$. We therefore choose
\[
 q(t)=-t^2\frac{\mathrm d}{\mathrm dt}\!\left(\frac{\Omega(t)}{t}\right)
 =\Omega(t)-t\Omega'(t).
\]
The large-displacement term then equals
$\Omega(r)-(r/R)\Omega(R)\le\Omega(r)$. At the same time, this choice
reduces the integral defining the slope by only $\Omega(1)-\Omega(r)$,
so the increments retain the required order. We prove that the essential
oscillation of $(-\Delta_p)^s v$ is comparable to $\omega(r)$ at every
sufficiently small scale $r$. Taking $u(x)=v(x_1)$ gives the example
in Theorem~\ref{thm:sharp-modulus}.

The paper is organized as follows. Section~\ref{sec:preliminaries}
introduces the operators and derives the viscosity inequalities used
for weak solutions. Section~\ref{sec:method} uses reflection to estimate
the integral on the annular cone and controls the localization and far-field terms.
Section~\ref{sec:estimates} constructs the concave profile $\phi$
and proves the interior estimates. Section~\ref{sec:sharpness} constructs
the sharpness examples and proves the Dini criterion.

\section{Preliminaries}\label{sec:preliminaries}

\subsection{Operators and elementary inequalities}\label{subsec:operators}

To prove the estimates in Theorems~\ref{thm:log} and~\ref{thm:logmodulus},
we work with translation-invariant operators of the form
\begin{equation}\label{eq:operator}
 \LL u(x)=\PV\int_{\R^n}J_p(u(x)-u(x+z))K(z)\,\mathrm dz,
 \qquad J_p(t)=|t|^{p-2}t.
\end{equation}
As in \cite[Section~1.2]{BT}, we assume that the kernel $K$ satisfies
\begin{equation}\label{eq:kernel}
 \begin{gathered}
 K\in C(\R^n\setminus\{0\}),\qquad K(z)=K(-z),\\
 \lambda |z|^{-n-sp}\le K(z)\le\Lambda |z|^{-n-sp}
 \quad(z\ne0),\qquad 0<\lambda\le\Lambda<\infty.
 \end{gathered}
\end{equation}
The choice $K(z)=|z|^{-n-sp}$ gives the fractional $p$-Laplacian
\eqref{eq:fracplap}.
Write $k(z)=|z|^{n+sp}K(z)$ and fix the near-field radius $\rho=1/16$.
For $t>0$, the oscillation of $k$ outside $B_\rho(0)$ is
\begin{equation}\label{eq:kernelmodulus}
 \omega_\rho(k;t)=\sup\bigl\{|k(z)-k(w)| : |z|,|w|\ge\rho,\ |z-w|\le t\bigr\}.
\end{equation}

We write $B_r=B_r(0)$.
The constants $c,C>0$ may change from line to line and depend only on
$n,p,\lambda,\Lambda$ and the additional parameters indicated.
We write $A\asymp B$ if $cB\le A\le CB$ for such constants.

The scalar identity below is given in \cite[Lemma~2.4]{BM}, and the
estimates follow from \cite[Lemmas~3.2--3.4]{KKL}.

\begin{lemma}\label{lem:scalar}
For $p>1$ and $(a,b)\ne(0,0)$,
\[
 J_p(a+b)-J_p(a)
 =(p-1)b\int_0^1|a+tb|^{p-2}\,\mathrm dt.
\]
The integral is comparable to $(|a|+|b|)^{p-2}$, with constants
depending only on $p$. Consequently, for all $a,b\in\R$,
\[
 |J_p(a)-J_p(b)|\le
 \begin{cases}
 C(|a|+|b|)^{p-2}|a-b|,&p\ge2,\\
 C|a-b|^{p-1},&1<p<2.
 \end{cases}
\]
\end{lemma}

The following elementary inequality gives a lower bound for $J_p(a)-J_p(b)$ in terms of $a-b$.

\begin{lemma}\label{lem:truncatedmonotonicity}
Let $p>1$ and $0<c_0\le C_0$. If $a,b\in\R$ satisfy $a\ge b$ and
\[
 [b,a]\cap\bigl([-C_0,-c_0]\cup[c_0,C_0]\bigr)\ne\varnothing,
\]
then
\[
 J_p(a)-J_p(b)\ge c\min\{a-b,1\},
\]
where $c>0$ depends only on $p,c_0,C_0$.
\end{lemma}
\begin{proof}
The case $a=b$ is immediate. For $a>b$, choose $t_0\in[b,a]$ with
$c_0\le|t_0|\le C_0$ and a subinterval $[t,t+h]\subset[b,a]$ containing
$t_0$, of length $h=\min\{a-b,1\}$. Since $c_0\le |t|+h\le C_0+2$,
monotonicity and Lemma~\ref{lem:scalar} give
\[
 J_p(a)-J_p(b)\ge J_p(t+h)-J_p(t)
 \asymp h(|t|+h)^{p-2}\ge c h. \qedhere
\]
\end{proof}

\subsection{Weak solutions}\label{subsec:weak}

\begin{definition}[Local weak solutions]\label{def:weak}
Let $U\subset\R^n$ be a domain and $f\in L^\infty_{\mathrm{loc}}(U)$.
A function $u\in W^{s,p}_{\mathrm{loc}}(U)\cap L^{p-1}_{sp}(\R^n)$
is a local weak subsolution (resp.\ supersolution) of $\LL u=f$ in $U$ if
\begin{equation}\label{eq:weakform}
 \frac12\iint_{\R^n\times\R^n}
 J_p(u(x)-u(y))(\zeta(x)-\zeta(y))K(y-x)
 \,\mathrm dx\,\mathrm dy
 \mathrel{\le\,(\ge)}\int_U f\zeta\,\mathrm dx
\end{equation}
for every nonnegative $\zeta\in C_c^\infty(U)$.
A local weak solution is both a subsolution and a supersolution.
Equivalently, \eqref{eq:weakform} holds with equality for every
$\zeta\in C_c^\infty(U)$.
\end{definition}

\begin{lemma}\label{lem:weakregularity}
Let $K$ satisfy \eqref{eq:kernel}. If
$u\in W^{s,p}_{\mathrm{loc}}(U)\cap L^{p-1}_{sp}(\R^n)$ is a local weak
solution of $\LL u=f$ in $U$ with $f\in L^\infty_{\mathrm{loc}}(U)$,
then $u$ is locally bounded and locally H\"older continuous in $U$.
\end{lemma}
\begin{proof}
If $sp>n$, the Sobolev--Morrey embedding \cite[Theorem~8.2]{DNPV} gives
$u\in C^{0,s-n/p}(\overline{B})$ for every ball $B\Subset U$.

Suppose $sp\le n$ and fix a ball $B\Subset U$. Choose
$\chi\in C_c^\infty(U)$ with
$0\le\chi\le1$ and $\chi=1$ near $\overline B$.
Then $v=\chi u$ belongs to $W^{s,p}(\R^n)$ and agrees with $u$ on $B$.
Only the interactions with points outside $\{\chi=1\}$ change in
\eqref{eq:weakform}. By symmetry of $K$, the resulting equation on
$B$ is $\LL v=f_\chi(x,v)$, where
\[
 f_\chi(x,t)=f(x)+\int_{\{\chi\ne1\}}
 \bigl[J_p(t-v(y))-J_p(t-u(y))\bigr]K(y-x)\,\mathrm dy.
\]
Since $\{\chi\ne1\}$ has positive distance from $B$, the kernel in
this integral is bounded by $C(1+|y|)^{-n-sp}$ uniformly for $x\in B$.
The tail assumption therefore gives, for a.e.\ $x\in B$ and all $t\in\R$,
\[
 \begin{aligned}
 |f_\chi(x,t)|&\le \norm{f}_{L^\infty(B)}
   +C\int_{\{\chi\ne1\}}
     \frac{|t|^{p-1}+|u(y)|^{p-1}}{(1+|y|)^{n+sp}}\,\mathrm dy\\
 &\le C(1+|t|^{p-1}).
 \end{aligned}
\]
Thus $f_\chi(\cdot,v)\in L^{p/(p-1)}(B)$, and by density the weak
identity extends to $W^{s,p}(\R^n)$ test functions compactly supported
in $B$. In the formulation of \cite[(2.10)]{Cozzi}, the right-hand side
is $-f_\chi$, and the preceding bound gives \cite[(2.11)]{Cozzi} with
$q=p$. It follows from \cite[Theorems~2.3 and~2.4]{Cozzi} that $v=u$
is locally bounded and locally H\"older continuous in $B$.
\end{proof}

\subsection{Viscosity solutions}
\label{subsec:viscosity}

We use the standard viscosity formulation (see \cite{BT,KKL}).

\begin{definition}[Viscosity solutions]\label{def:viscosity}
Let $U\subset\R^n$ be a domain, $f:U\to\R$ locally bounded, and
$u:\R^n\to\R$ measurable. For a ball $B_r(x)\Subset U$ and a function
$\varphi$ of class $C^2$ on a neighborhood of $\overline{B}_r(x)$, set
\[
 v(w)=
 \begin{cases}
 \varphi(w),&w\in B_r(x),\\
 u(w),&w\notin B_r(x).
 \end{cases}
\]
We call $u$ a viscosity subsolution (resp.\ supersolution) of $\LL u=f$
in $U$ if $u$ is upper (resp.\ lower) semicontinuous in $U$,
$u_+\in L^{p-1}_{sp}(\R^n)$ (resp.\ $u_-\in L^{p-1}_{sp}(\R^n)$),
and for every such $B_r(x)$ and $\varphi$ with $\varphi(x)=u(x)$ and
$\varphi(w)\mathrel{\ge\,(\le)}u(w)$ for all $w\in B_r(x)$, we have
\[
 \LL v(x)\mathrel{\le\,(\ge)}f(x).
\]
A viscosity solution is both a subsolution and a supersolution.
\end{definition}

If $\nabla\varphi(x)=0$, the local integral converges absolutely since
$2(p-1)>sp$. For test functions with nonzero gradient we use the
following principal value estimate.

\begin{lemma}\label{lem:smalltest}
Let $K$ satisfy \eqref{eq:kernel}, and let $v\in C^2(B_r(x))$ have bounded
Hessian and $\nabla v(x)\ne0$. There is a constant $C=C(n,p)$ such that
for every measurable set $E\subset B_r$ with $E=-E$, the principal value
exists and satisfies
\begin{equation}\label{eq:testbound}
 \begin{aligned}
 &\left|\PV\int_E J_p(v(x)-v(x+z))K(z)\,\mathrm dz\right|\\
 &\qquad\le C\Lambda r\,\norm{D^2v}_{L^\infty(B_r(x))}
       \bigl(|\nabla v(x)|+r\,\norm{D^2v}_{L^\infty(B_r(x))}\bigr)^{p-2}.
 \end{aligned}
\end{equation}
\end{lemma}
\begin{proof}
Subtract $J_p(-\nabla v(x)\cdot z)K(z)$, whose integral over
$E\setminus B_\varepsilon$ vanishes by oddness. Spherical integration
as in \cite[Lemma~3.5 and (3.2)]{KKL}, with $p-sp=1$, bounds the
$L^1(B_r)$ norm of the remaining integrand by
\[
 C\Lambda\norm{D^2v}_{L^\infty(B_r(x))}
 \int_0^r\bigl(|\nabla v(x)|+t\norm{D^2v}_{L^\infty(B_r(x))}\bigr)^{p-2}
 \,\mathrm dt.
\]
The remainder is therefore integrable on $E$. Rescaling the radial
integral to $(0,1)$ and applying Lemma~\ref{lem:scalar} gives
\eqref{eq:testbound}.
\end{proof}

\begin{proposition}\label{prop:weakviscosity}
Let $K$ satisfy \eqref{eq:kernel}, and let
$u\in W^{s,p}_{\mathrm{loc}}(U)\cap L^{p-1}_{sp}(\R^n)$ be a
local weak solution of $\LL u=f$ in $U$, with
$f\in L^\infty_{\mathrm{loc}}(U)$.
On every ball $B\Subset U$, $u$ satisfies
\[
 \operatorname*{ess\,inf}_B f\le\LL u\le\operatorname*{ess\,sup}_B f
\]
in the viscosity sense. If $f$ is continuous, $u$ is a viscosity solution
of $\LL u=f$ in $U$.
\end{proposition}
\begin{proof}
Fix a ball $B\Subset U$. By Lemma~\ref{lem:weakregularity},
$u$ is continuous on $\overline B$. Testing \eqref{eq:weakform} with
nonnegative functions in $C_c^\infty(B)$ shows that $u$ is a weak
subsolution with datum $\operatorname*{ess\,sup}_B f$ and a weak
supersolution with datum $\operatorname*{ess\,inf}_B f$ on $B$.
Both data are constant, so Proposition~1.5 of \cite{BT} gives the
corresponding viscosity inequalities on $B$. If $f$ is continuous,
the same proposition applies directly to $\LL u=f$ on each such ball
(see also \cite[Theorem~1.3]{BM} for the fractional $p$-Laplacian
when $n>sp$).
\end{proof}

\section{Integral estimates at a positive maximum}\label{sec:method}

Let $\phi$ be a candidate modulus of continuity for $u$ and let $L>0$.
To keep the first contact point in the interior, we use a nonnegative
function $\psi$ that vanishes on $B_1$ and increases with $|x|$ outside $B_1$.
It suffices to rule out a positive maximum of $u(x)-u(y)-L\phi(|x-y|)-\psi(x)$
for $x,y\in\overline{B}_{11/8}$. Restricting to $B_1$ and interchanging $x$
and $y$ then gives $|u(x)-u(y)|\le L\phi(|x-y|)$.
In this section we leave $\phi$ unspecified and estimate the near- and
far-field integrals at a positive maximum. In Section~\ref{sec:estimates}
we choose $\phi$ according to the oscillation bound and use these estimates
to rule out such a maximum.

Suppose that $u\in L^\infty(B_{3/2})\cap L^{p-1}_{sp}(\R^n)$ satisfies
$-F\le\LL u\le F$ in $B_{3/2}$ in the viscosity sense for some $F\ge0$; by
Definition~\ref{def:viscosity}, $u$ is continuous in $B_{3/2}$. Set
\[
 A=\norm{u}_{L^\infty(B_{3/2})}+\Tail(u;0,3/2)+F^{1/(p-1)}.
\]
If $A=0$, the estimates are immediate. For $A>0$, homogeneity gives
$\LL(u/A)=A^{1-p}\LL u$ in the viscosity sense. Thus, when $\LL u=f$,
the right-hand side and its modulus both scale by $A^{1-p}$.
Writing $u$ for $u/A$, we have
\begin{equation}\label{eq:normalized}
 \norm{u}_{L^\infty(B_{3/2})}+\Tail(u;0,3/2)\le1,
 \qquad -1\le\LL u\le1\quad\text{in }B_{3/2},
\end{equation}
where the inequalities for $\LL u$ hold in the viscosity sense.

\subsection{The contact inequality at a positive maximum}

Let $\phi$ be continuous and nondecreasing on $[0,3]$, with $\phi(0)=0$.
For $x,y\in\overline{B}_{11/8}$, we have $|x-y|\le11/4<3$, so
$\phi(|x-y|)$ is well defined.
On $(0,\rho)$, where $\rho=1/16$, assume that $\phi$ is $C^2$ with
$\inf_{(0,\rho)}\phi'>0$.
Set
\begin{equation}\label{eq:localization}
 \psi(x)=m[(|x|^2-1)_+]^3,\qquad
 \Phi(x,y)=L\phi(|x-y|)+\psi(x).
\end{equation}
Here $m>0$ is a fixed numerical constant chosen so that $\psi(x)>2$ for
$|x|\ge5/4$. The cubic power makes $\psi$ of class $C^2$ across $\partial B_1$.
For the remainder of this section, suppose the maximum is positive and
denote by $(\bar x,\bar y)$ any pair attaining it, so that
\begin{equation}\label{eq:max}
 \max_{x,y\in\overline{B}_{11/8}}\{u(x)-u(y)-\Phi(x,y)\}
 =u(\bar x)-u(\bar y)-\Phi(\bar x,\bar y)>0.
\end{equation}
By \eqref{eq:max} and \eqref{eq:normalized},
$\Phi(\bar x,\bar y)<u(\bar x)-u(\bar y)\le2$.
The maximizing points are distinct, since $u(x)-u(x)-\Phi(x,x)=-\psi(x)\le0$.
Moreover, $\psi(\bar x)<2$ forces $|\bar x|<5/4$.
For the rest of this section, take $L$ large enough that $L\phi(8\rho/9)>2$
and $L\inf_{(0,\rho)}\phi'\ge2\|\nabla\psi\|_{L^\infty(B_{11/8})}$.
Since $\phi$ is nondecreasing,
a positive maximum then requires $|\bar x-\bar y|<8\rho/9$.
Altogether, the maximizing pair satisfies
\begin{equation}\label{eq:maxproperties}
 \begin{gathered}
 L\phi(|\bar x-\bar y|)+\psi(\bar x)<u(\bar x)-u(\bar y)\le2,\\
 |\bar x|<5/4,\qquad 0<|\bar x-\bar y|<8\rho/9.
 \end{gathered}
\end{equation}
In particular, $|\bar x-\bar y|\to0$ uniformly as $L\to\infty$.

Since $|\bar y|<5/4+\rho$ and $5/4+2\rho=11/8$, both $B_\rho(\bar x)$
and $B_\rho(\bar y)$ lie in $B_{11/8}$. At each positive maximum $(\bar x,\bar y)$,
we take the common radius of the two test balls to be
\begin{equation}\label{eq:delta}
 \delta=|\bar x-\bar y|^2.
\end{equation}
By \eqref{eq:maxproperties}, $0<\delta<|\bar x-\bar y|/8$.
As in Definition~\ref{def:viscosity}, define the globally patched test functions
\[
 v_{\bar x}(w):=
 \begin{cases}
 u(\bar x)+\Phi(w,\bar y)-\Phi(\bar x,\bar y),&w\in B_\delta(\bar x),\\
 u(w),&w\notin B_\delta(\bar x),
 \end{cases}
\]
and
\[
 v_{\bar y}(w):=
 \begin{cases}
 u(\bar y)-\Phi(\bar x,w)+\Phi(\bar x,\bar y),&w\in B_\delta(\bar y),\\
 u(w),&w\notin B_\delta(\bar y).
 \end{cases}
\]
By \eqref{eq:max}, $v_{\bar x}\ge u$ on $B_\delta(\bar x)$ and
$v_{\bar y}\le u$ on $B_\delta(\bar y)$, with equality at the respective centers.
Since $\delta<|\bar x-\bar y|/8$ and $9|\bar x-\bar y|/8<\rho$,
the smooth branches are $C^2$ on neighborhoods of the respective closed balls.

At the contact points, $|\nabla v_{\bar y}(\bar y)|=L\phi'(|\bar x-\bar y|)$
and $\nabla v_{\bar x}(\bar x)=\nabla v_{\bar y}(\bar y)+\nabla\psi(\bar x)$.
Since $|\nabla\psi(\bar x)|\le\tfrac12 L\phi'(|\bar x-\bar y|)$, we have
\[
 \tfrac12 L\phi'(|\bar x-\bar y|)\le|\nabla v_{\bar x}(\bar x)|
 \le\tfrac32 L\phi'(|\bar x-\bar y|).
\]
Both gradients are nonzero, so Lemma~\ref{lem:smalltest} gives principal
value convergence at both points. Applying \eqref{eq:normalized} to the
upper test at $\bar x$ and the lower test at $\bar y$ gives
$\LL v_{\bar x}(\bar x)\le1$ and $\LL v_{\bar y}(\bar y)\ge-1$.
Hence $\LL v_{\bar x}(\bar x)-\LL v_{\bar y}(\bar y)\le2$.

For $|z|<\delta$, the increments are given by the smooth tests.
For $\delta\le |z|<\rho$, the increments are those of $u$ and both translated
points remain in $B_{11/8}$, so \eqref{eq:max} applies. For $|z|\ge\rho$,
the translated points may leave $B_{11/8}$.

Set $\Delta_z w(x):=w(x)-w(x+z)$.
For a measurable set $V\subset\R^n$ such that
$V\cap B_r=-(V\cap B_r)$ for some $r>0$, the contribution of
displacements $z\in V$ to the operator difference is
\[
 I_V=\PV\int_V\bigl[J_p(\Delta_z v_{\bar x}(\bar x))
                     -J_p(\Delta_z v_{\bar y}(\bar y))\bigr]K(z)\,\mathrm dz.
\]
Write $I_{\mathrm{far}}=I_{\R^n\setminus B_\rho}$. Splitting the integral over the three regions gives
\begin{equation}\label{eq:differenceineq}
 \LL v_{\bar x}(\bar x)-\LL v_{\bar y}(\bar y)
 = I_{B_\delta}+I_{B_\rho\setminus B_\delta}+I_{\mathrm{far}}\le2.
\end{equation}
If $u$ is a viscosity solution of $\LL u=f$, the same two contact inequalities
give the upper bound $f(\bar x)-f(\bar y)$ in \eqref{eq:differenceineq}.

\subsection{Estimates for the integral terms}

In \eqref{eq:differenceineq}, we seek a positive contribution to
$I_{B_\rho\setminus B_\delta}$ from displacements directed along the line
through $\bar x$ and $\bar y$. This contribution must dominate the errors
in the annulus and any negative contributions from $I_{B_\delta}$ and
$I_{\mathrm{far}}$, so that the sum exceeds the upper bound in
\eqref{eq:differenceineq}. For the remainder of this section, assume that
$\phi$ is concave on $(0,\rho)$. To quantify its concavity at the contact scale
$|\bar x-\bar y|$, we use
\[
 \gamma_\phi(t):=-\frac{t\phi''(t)}{\phi'(t)}
                =-\frac{\mathrm d\log\phi'(t)}{\mathrm d\log t},
 \qquad 0<t<\rho.
\]
We also assume that there exists $M_\phi\ge1$ such that
\begin{equation}\label{eq:profileconditions}
 \begin{aligned}
 \phi'(r)\ge M_\phi^{-1},\quad 0<\gamma_\phi(r)\le M_\phi, &\qquad r\in(0,\rho),\\
 M_\phi^{-1}\gamma_\phi(r)\le \gamma_\phi(t)\le M_\phi \gamma_\phi(r), &\qquad
 r\in(0,8\rho/9),\quad t\in[7r/8,9r/8].
 \end{aligned}
\end{equation}

Throughout this subsection, the threshold for $L$ depends only on $n,p,\lambda,\Lambda$, and $M_\phi$.

On $B_\rho\setminus B_\delta$, the integral in \eqref{eq:differenceineq} involves the
increments of $u$ at both contact points. In \eqref{eq:max}, the
localization term $\psi$ depends only on the first variable. Translating both
points by the same $z\in B_\rho$ leaves their distance unchanged, so maximality gives
\[
 \Delta_z(u-\psi)(\bar x)
 =\Delta_z u(\bar x)-\Delta_z\psi(\bar x)
 \ge\Delta_z u(\bar y).
\]
We use this inequality to obtain the lower bound on the cone in
Lemma~\ref{lem:fourpoint}. Lemma~\ref{lem:localerror} estimates the error from
replacing $J_p(\Delta_z u(\bar x))$ by $J_p(\Delta_z(u-\psi)(\bar x))$.

To estimate the annular term $I_{B_\rho\setminus B_\delta}$ in
\eqref{eq:differenceineq} from below, we first consider displacements in a cone
about the line through $\bar x$ and $\bar y$. Set
$\bar e=(\bar x-\bar y)/|\bar x-\bar y|$ and $\eta=1/(64M_\phi)$, and let
\[
 E:=\bigl\{z\in\R^n:\ \tfrac12\eta|\bar x-\bar y|<|z|<\eta|\bar x-\bar y|,\;
                         |\bar e\cdot z|\ge |z|/2\bigr\}.
\]

\begin{lemma}\label{lem:fourpoint}
Assume that $u\in C(B_{3/2})\cap L^{p-1}_{sp}(\R^n)$ satisfies \eqref{eq:normalized}
and $\phi$ satisfies \eqref{eq:profileconditions}. For all sufficiently large $L$
and every pair $(\bar x,\bar y)$ satisfying \eqref{eq:max}, we have
\begin{equation}\label{eq:roughgain}
 \int_E \bigl[J_p(\Delta_z(u-\psi)(\bar x))-J_p(\Delta_z u(\bar y))\bigr]K(z)\,\mathrm dz
 \ge c\bigl(L\phi'(|\bar x-\bar y|)\bigr)^{p-1}\gamma_\phi(|\bar x-\bar y|).
\end{equation}
Here $c>0$ depends only on $n,p,\lambda,M_\phi$.
\end{lemma}
\begin{proof}
Set $d(z)=\Delta_z(u-\psi)(\bar x)-\Delta_z u(\bar y)\ge0$ for $z\in B_\rho$.
Let $E_+=\{z\in E:\bar e\cdot z>0\}$
and set $Rz=z-2(\bar e\cdot z)\bar e$. Then $E$ is the disjoint union of
$E_+$ and $R(E_+)$. For $z\in E_+$ set $t=\bar e\cdot z>0$.
Since $z-Rz=2t\bar e$ and $\bar x-\bar y=|\bar x-\bar y|\bar e$, we have
\[
 |(\bar x+z)-(\bar y+Rz)|=|\bar x-\bar y|+2t,\qquad |(\bar x+Rz)-(\bar y+z)|=|\bar x-\bar y|-2t.
\]
Comparing these two crossed pairs with $(\bar x,\bar y)$ in \eqref{eq:max} gives
\[
 \begin{aligned}
 \Delta_z(u-\psi)(\bar x)-\Delta_{Rz}u(\bar y)&\ge L[\phi(|\bar x-\bar y|)-\phi(|\bar x-\bar y|+2t)],\\
 \Delta_{Rz}(u-\psi)(\bar x)-\Delta_z u(\bar y)&\ge L[\phi(|\bar x-\bar y|)-\phi(|\bar x-\bar y|-2t)].
 \end{aligned}
\]
Adding these inequalities gives
\begin{equation}\label{eq:crosssum}
 \begin{aligned}
 d(z)+d(Rz)
 &\ge L\bigl[2\phi(|\bar x-\bar y|)-\phi(|\bar x-\bar y|+2t)-\phi(|\bar x-\bar y|-2t)\bigr]\\
 &\ge cL\phi'(|\bar x-\bar y|)\gamma_\phi(|\bar x-\bar y|)\frac{|z|^2}{|\bar x-\bar y|}.
 \end{aligned}
\end{equation}
Indeed, $9|\bar x-\bar y|/8<\rho$ by \eqref{eq:maxproperties} and $2t<|\bar x-\bar y|/8$.
Writing the second difference as
\[
 \int_0^{2t}(2t-h)\bigl[-\phi''(|\bar x-\bar y|+h)-\phi''(|\bar x-\bar y|-h)\bigr]\,\mathrm dh
\]
and applying \eqref{eq:profileconditions} at the scales $|\bar x-\bar y|$
and $|\bar x-\bar y|\pm h$ gives the second inequality in \eqref{eq:crosssum},
since $t\ge |z|/2$.

Comparing the maximum at $(\bar x,\bar y)$ with the
values at $(\bar x+z,\bar y)$ and $(\bar x,\bar y+z)$ gives
\[
 \begin{aligned}
 \Delta_z(u-\psi)(\bar x)&\ge L[\phi(|\bar x-\bar y|)-\phi(|\bar x-\bar y+z|)]
 \ge -CL\phi'(|\bar x-\bar y|)|z|,\\
 \Delta_z u(\bar y)&\le L[\phi(|\bar x-\bar y-z|)-\phi(|\bar x-\bar y|)]
 \le -cL\phi'(|\bar x-\bar y|)|z|.
 \end{aligned}
\]
Indeed,
\[
 |\bar x-\bar y|-|\bar x-\bar y-z|=\frac{2|\bar x-\bar y|t-|z|^2}{|\bar x-\bar y|+|\bar x-\bar y-z|}
 \ge\frac{1-\eta}{2+\eta}|z|.
\]
Since
$|\bar x-\bar y+Rz|=|\bar x-\bar y-z|$ and $|\bar x-\bar y-Rz|=|\bar x-\bar y+z|$, the reflected pair satisfies
\[
 \Delta_{Rz}(u-\psi)(\bar x)\ge cL\phi'(|\bar x-\bar y|)|z|,\qquad
 \Delta_{Rz}u(\bar y)\le CL\phi'(|\bar x-\bar y|)|z|.
\]

Applying $J_p$ to the two increments gives the difference
\[
 d_p(z):=J_p\bigl(\Delta_z(u-\psi)(\bar x)\bigr)
         -J_p\bigl(\Delta_z u(\bar y)\bigr).
\]
Since $d(z)\ge0$ and $J_p$ is increasing, $d_p(z)\ge0$ on $B_\rho$.
For $z\in E$, normalize the increments by setting
\[
 a_z=\frac{\Delta_z(u-\psi)(\bar x)}{L\phi'(|\bar x-\bar y|)|z|},
 \qquad
 b_z=\frac{\Delta_z u(\bar y)}{L\phi'(|\bar x-\bar y|)|z|}.
\]
For $z\in E_+$, the increment bounds show that $[b_z,a_z]$ meets $[-C_0,-c_0]$ for some
$0<c_0\le C_0$ depending only on $M_\phi$. For $Rz$ the same argument
gives an intersection with $[c_0,C_0]$.
Lemma~\ref{lem:truncatedmonotonicity} therefore gives
\[
 \begin{aligned}
 J_p(a_z)-J_p(b_z)&\ge c\min\{a_z-b_z,1\},\\
 J_p(a_{Rz})-J_p(b_{Rz})&\ge c\min\{a_{Rz}-b_{Rz},1\}.
 \end{aligned}
\]
Since $|Rz|=|z|$, the two increments have the same normalization.
Adding the estimates, using $\min\{a,1\}+\min\{b,1\}\ge\min\{a+b,1\}$
for $a,b\ge0$, and rescaling by the $(p-1)$-homogeneity of $J_p$, we obtain
\[
 \begin{aligned}
 d_p(z)+d_p(Rz)
 &\ge c\bigl(L\phi'(|\bar x-\bar y|)|z|\bigr)^{p-2}\\
 &\qquad{}\times\min\bigl\{d(z)+d(Rz),\,L\phi'(|\bar x-\bar y|)|z|\bigr\}.
 \end{aligned}
\]
On $E$ we have $\gamma_\phi(|\bar x-\bar y|)|z|/|\bar x-\bar y|\le M_\phi\eta\le1$.
By \eqref{eq:crosssum} the last minimum is at least
\[
 cL\phi'(|\bar x-\bar y|)\gamma_\phi(|\bar x-\bar y|)
 \frac{|z|^2}{|\bar x-\bar y|}.
\]
Consequently, for $z\in E_+$,
\begin{equation}\label{eq:fourpointgain}
 d_p(z)+d_p(Rz)
 \ge c\bigl(L\phi'(|\bar x-\bar y|)\bigr)^{p-1}\gamma_\phi(|\bar x-\bar y|)\frac{|z|^p}{|\bar x-\bar y|}.
\end{equation}

Since $d_p\ge0$, we may use $K(z)\ge\lambda|z|^{-n-sp}$
before reflecting in the integral. This gives
\[
 \int_E d_p(z)K(z)\,\mathrm dz
 \ge\lambda\int_{E_+}\frac{d_p(z)+d_p(Rz)}
                              {|z|^{n+sp}}\,\mathrm dz .
\]
The angular measure of $E_+$ is fixed. Since $p-sp=1$, the radial factor
in the integral of \eqref{eq:fourpointgain} is
\[
 \frac{1}{|\bar x-\bar y|}\int_{\eta|\bar x-\bar y|/2}^{\eta|\bar x-\bar y|}r^{p-sp-1}\,\mathrm dr
 =\frac{\eta}{2}.
\]
This proves \eqref{eq:roughgain}.
\end{proof}

\begin{lemma}\label{lem:localerror}
Assume that $u\in C(B_{3/2})\cap L^{p-1}_{sp}(\R^n)$ satisfies \eqref{eq:normalized}
and $\phi$ satisfies \eqref{eq:profileconditions}. For all sufficiently large $L$
and every pair $(\bar x,\bar y)$ satisfying \eqref{eq:max}, we have
\begin{equation}\label{eq:annulus}
 I_{B_\rho\setminus B_\delta}
 \ge \int_E\bigl[J_p(\Delta_z(u-\psi)(\bar x))-J_p(\Delta_z u(\bar y))\bigr]
 K(z)\,\mathrm dz-C|\bar x-\bar y|^{\sigma}.
\end{equation}
Here $\sigma=\sigma(p)>0$ and $C$ depends only on $n,p,\lambda,\Lambda,M_\phi$.
\end{lemma}
\begin{proof}
Let $\delta<r<\rho$ be such that $E\subset B_r\setminus B_\delta$.
By \eqref{eq:max}, $\Delta_z(u-\psi)(\bar x)\ge\Delta_z u(\bar y)$ for $z\in B_\rho$.
Since $J_p$ is increasing, adding and subtracting $J_p(\Delta_z(u-\psi)(\bar x))$
on $B_r\setminus B_\delta$ and restricting the nonnegative term to $E$ gives
\begin{equation}\label{eq:annularsplit}
 \begin{aligned}
 I_{B_\rho\setminus B_\delta}
 &\ge \int_E\bigl[J_p(\Delta_z(u-\psi)(\bar x))-J_p(\Delta_z u(\bar y))\bigr]
 K(z)\,\mathrm dz\\
 &\quad-\int_{B_r\setminus B_\delta}
 \bigl|J_p(\Delta_z u(\bar x))-J_p(\Delta_z(u-\psi)(\bar x))\bigr|
 K(z)\,\mathrm dz\\
 &\quad-|I_{B_\rho\setminus B_r}|.
 \end{aligned}
\end{equation}
It suffices to show that for some such $r$ both error terms in \eqref{eq:annularsplit} are bounded by
$C|\bar x-\bar y|^\sigma$ for some $\sigma=\sigma(p)>0$.

Fix $0<\kappa<1$. For $x_0\in B_{11/8}$ we have
$B_{2\rho}(x_0)\Subset B_{3/2}$. Applying \cite[Theorem~2.1 and
Remark~2.1(2)]{BT} to $u(x_0+\rho\,\cdot)$ gives
\[
 [u]_{C^{0,\kappa}(B_\rho(x_0))}
 \le \frac{C}{\rho^\kappa}
 \left(\norm{u}_{L^\infty(B_{2\rho}(x_0))}
       +\Tail(u;x_0,2\rho)+\rho\right),
\]
where we used $-1\le\LL u\le1$ from \eqref{eq:normalized}.
After rescaling, the right-hand side is bounded by $\rho^{sp}$, giving
the term $(\rho^{sp})^{1/(p-1)}=\rho$ since $sp=p-1$.
The norm and tail are uniformly bounded because $|u|\le1$ on $B_{3/2}$,
$|w-x_0|\ge|w|/12$ for $|w|\ge3/2$, and $\Tail(u;0,3/2)\le1$.
Since $|u(x)-u(y)|\le2$ when $|x-y|\ge\rho$, we obtain
\begin{equation}\label{eq:preholder}
 [u]_{C^{0,\kappa}(B_{11/8})}\le C,
\end{equation}
where $C$ may also depend on $\kappa$.

We now estimate $|I_{B_\rho\setminus B_r}|$. For $|z|<\rho$, applying
\eqref{eq:preholder} to the pairs $(\bar x,\bar y)$ and $(\bar x+z,\bar y+z)$ gives
\[
 |\Delta_z u(\bar x)-\Delta_z u(\bar y)|
 \le |u(\bar x)-u(\bar y)|+|u(\bar x+z)-u(\bar y+z)|
 \le C|\bar x-\bar y|^\kappa.
\]
Lemma~\ref{lem:scalar} and $|u|\le1$ then give
\[
 |J_p(\Delta_z u(\bar x))-J_p(\Delta_z u(\bar y))|
 \le C|\bar x-\bar y|^{\kappa\min\{1,p-1\}}.
\]
As $\int_{|z|\ge r}K(z)\,\mathrm dz\le Cr^{-sp}$, we retain half
this power by taking
\[
 r=|\bar x-\bar y|^\theta,\qquad
 \theta=\frac{\kappa\min\{1,p-1\}}{2sp}.
\]
Indeed, $sp=p-1$ gives $0<\theta<1$, so \eqref{eq:maxproperties} and
\eqref{eq:delta} yield
$\delta<\tfrac12\eta|\bar x-\bar y|<\eta|\bar x-\bar y|<r<\rho$
for all sufficiently large $L$. Thus the inner annulus contains $E$, and
\begin{equation}\label{eq:outerannulus}
 |I_{B_\rho\setminus B_r}|
 \le C|\bar x-\bar y|^{\kappa\min\{1,p-1\}}r^{-sp}
 =C|\bar x-\bar y|^{\kappa\min\{1,p-1\}/2}.
\end{equation}

For the localization error on $B_r\setminus B_\delta$, we now take $\kappa$
sufficiently close to $1$, depending only on $p$, so that $(1-\kappa)(p-2)_+<\kappa/3$.
By \eqref{eq:maxproperties} and \eqref{eq:preholder},
$\psi(\bar x)\le u(\bar x)-u(\bar y)\le C|\bar x-\bar y|^\kappa$.
By \eqref{eq:localization}, $\nabla\psi(x)=6m^{1/3}x\,\psi(x)^{2/3}$.
Using this identity at $\bar x$ and Taylor's formula with $D^2\psi$ bounded on $B_{11/8}$, we obtain
\begin{equation}\label{eq:incrementbounds}
 \begin{aligned}
 |\nabla\psi(\bar x)|&\le C|\bar x-\bar y|^{2\kappa/3},\\
 |\Delta_z\psi(\bar x)|&\le C\bigl(|z|^2+|\nabla\psi(\bar x)|\,|z|\bigr),
 \qquad |z|<r.
 \end{aligned}
\end{equation}
In the following estimates the increments are evaluated at $\bar x$.

For $p\ge2$, \eqref{eq:preholder} and $|\Delta_z\psi|\le C|z|$ give
$|\Delta_z u|+|\Delta_z(u-\psi)|\le C|z|^\kappa$. By Lemma~\ref{lem:scalar},
\[
 \begin{aligned}
 |J_p(\Delta_z u)-J_p(\Delta_z(u-\psi))|
 &\le C\bigl(|\Delta_z u|+|\Delta_z(u-\psi)|\bigr)^{p-2}|\Delta_z\psi|\\
 &\le C|z|^{\kappa(p-2)}\bigl(|z|^2+|\nabla\psi(\bar x)|\,|z|\bigr).
 \end{aligned}
\]
Put $\beta=(1-\kappa)(p-2)<\kappa/3$. Polar integration with
$K(z)\le\Lambda|z|^{-n-sp}$ and $sp=p-1$ gives
\[
 \begin{aligned}
 \int_{B_r\setminus B_\delta}
 |J_p(\Delta_z u)-J_p(\Delta_z(u-\psi))|K(z)\,\mathrm dz
 &\le C\int_\delta^rt^{-\beta}
 \left(1+\frac{|\nabla\psi(\bar x)|}{t}\right)\,\mathrm dt\\
 &\le Cr^{1-\beta}
 +C|\nabla\psi(\bar x)|\,\delta^{-\beta}\log(r/\delta).
 \end{aligned}
\]
Here $1-\beta>0$, and \eqref{eq:incrementbounds} and
$\delta=|\bar x-\bar y|^2$ give
\[
 |\nabla\psi(\bar x)|\,\delta^{-\beta}
 \le C|\bar x-\bar y|^{2\kappa/3-2\beta},\qquad 2\kappa/3-2\beta>0.
\]

For $1<p<2$, Lemma~\ref{lem:scalar} and \eqref{eq:incrementbounds} give
\[
 \begin{aligned}
 |J_p(\Delta_z u)-J_p(\Delta_z(u-\psi))|
 &\le C|\Delta_z\psi|^{p-1}\\
 &\le C\bigl(|z|^{2(p-1)}+|\nabla\psi(\bar x)|^{p-1}|z|^{p-1}\bigr).
 \end{aligned}
\]
Polar integration now yields
\[
 \begin{aligned}
 \int_{B_r\setminus B_\delta}
 |J_p(\Delta_z u)-J_p(\Delta_z(u-\psi))|K(z)\,\mathrm dz
 &\le C\int_\delta^r
 \left(t^{p-2}+\frac{|\nabla\psi(\bar x)|^{p-1}}{t}\right)\,\mathrm dt\\
 &\le Cr^{p-1}
 +C|\bar x-\bar y|^{2\kappa(p-1)/3}\log(r/\delta).
 \end{aligned}
\]
Since $r=|\bar x-\bar y|^\theta$ and
$\log(r/\delta)=(2-\theta)\log(1/|\bar x-\bar y|)$,
each power of $|\bar x-\bar y|$ in these bounds is positive and depends
only on $p$. Taking $\sigma=\sigma(p)>0$ smaller than these powers and
the exponent in \eqref{eq:outerannulus} absorbs the logarithm.
Substitution in \eqref{eq:annularsplit} proves \eqref{eq:annulus}.
\end{proof}

\begin{lemma}\label{lem:contactballs}
Assume that $K$ satisfies \eqref{eq:kernel} and that $\phi$ satisfies the slope
and concavity bounds in the first line of \eqref{eq:profileconditions}.
For all sufficiently large $L$ and every pair $(\bar x,\bar y)$ satisfying \eqref{eq:max}, we have
\begin{equation}\label{eq:patchbound}
 |I_{B_\delta}|\le C\bigl(L\phi'(|\bar x-\bar y|)\bigr)^{p-1}|\bar x-\bar y|.
\end{equation}
The constant $C$ and the threshold for $L$ depend only on $n,p,\Lambda,M_\phi$.
\end{lemma}
\begin{proof}
The argument of $\phi$ is $t=|x-\bar y|$ on $B_\delta(\bar x)$ and
$t=|\bar x-y|$ on $B_\delta(\bar y)$. In both cases $t$ differs from
$|\bar x-\bar y|$ by at most $\delta<|\bar x-\bar y|/8$.
Since \eqref{eq:profileconditions} gives $\lvert(\log\phi')'(t)\rvert\le M_\phi/t$, we have
\[
 \phi'(t)\asymp\phi'(|\bar x-\bar y|),\qquad
 |\phi''(t)|+\frac{\phi'(t)}{t}\le C\frac{\phi'(|\bar x-\bar y|)}{|\bar x-\bar y|}.
\]
Differentiating the two tests, using the boundedness of $D^2\psi$ and the gradient
bounds preceding \eqref{eq:differenceineq}, we obtain, for $\xi\in\{\bar x,\bar y\}$,
\[
 |\nabla v_\xi(\xi)|\asymp L\phi'(|\bar x-\bar y|),\qquad
 \norm{D^2v_\xi}_{L^\infty(B_\delta(\xi))}
 \le C\frac{L\phi'(|\bar x-\bar y|)}{|\bar x-\bar y|}.
\]
Thus $\delta\norm{D^2v_\xi}_{L^\infty(B_\delta(\xi))}\le C|\nabla v_\xi(\xi)|$.
Applying \eqref{eq:testbound} at both contact points and using
$\delta=|\bar x-\bar y|^2$ gives
\[
 |I_{B_\delta}|
 \le C\bigl(L\phi'(|\bar x-\bar y|)\bigr)^{p-1}\frac{\delta}{|\bar x-\bar y|}
 = C\bigl(L\phi'(|\bar x-\bar y|)\bigr)^{p-1}|\bar x-\bar y|. \qedhere
\]
\end{proof}

\begin{lemma}\label{lem:nearfield}
Assume that $u\in C(B_{3/2})\cap L^{p-1}_{sp}(\R^n)$ satisfies \eqref{eq:normalized}
and $\phi$ satisfies \eqref{eq:profileconditions}. For all sufficiently large $L$
and every pair $(\bar x,\bar y)$ satisfying \eqref{eq:max}, we have
\begin{equation}\label{eq:vanishing}
 I_{B_\rho}
 \ge c\bigl(L\phi'(|\bar x-\bar y|)\bigr)^{p-1}\gamma_\phi(|\bar x-\bar y|)
       -C|\bar x-\bar y|^{\sigma}.
\end{equation}
Here $\sigma=\sigma(p)>0$. The constant $c>0$ depends only on
$n,p,\lambda,M_\phi$, and $C$ depends only on $n,p,\lambda,\Lambda,M_\phi$.
\end{lemma}
\begin{proof}
Choose $0<\kappa<1$ depending only on $p$ so that $(1-\kappa)(p-1)\le1/2$.
Since $\phi$ is concave on $[0,\rho)$ and $\phi(0)=0$, we have $t\phi'(t)\le\phi(t)$
for $0<t<\rho$. By \eqref{eq:maxproperties}, $\psi\ge0$, and \eqref{eq:preholder},
\[
 L|\bar x-\bar y|\phi'(|\bar x-\bar y|)
 \le L\phi(|\bar x-\bar y|)<u(\bar x)-u(\bar y)\le C|\bar x-\bar y|^\kappa.
\]
Lemma~\ref{lem:contactballs} therefore gives
\[
 |I_{B_\delta}|
 \le C\bigl(L\phi'(|\bar x-\bar y|)\bigr)^{p-1}|\bar x-\bar y|
 \le C|\bar x-\bar y|^{1-(1-\kappa)(p-1)}
 \le C|\bar x-\bar y|^{1/2}.
\]
Taking $\sigma$ to be the minimum of $1/2$ and the exponent in Lemma~\ref{lem:localerror},
we obtain from Lemmas~\ref{lem:fourpoint} and~\ref{lem:localerror}
\begin{align*}
 I_{B_\rho}
 &=I_{B_\rho\setminus B_\delta}+I_{B_\delta}\\
 &\ge c\bigl(L\phi'(|\bar x-\bar y|)\bigr)^{p-1}\gamma_\phi(|\bar x-\bar y|)
      -C|\bar x-\bar y|^\sigma-|I_{B_\delta}|\\
 &\ge c\bigl(L\phi'(|\bar x-\bar y|)\bigr)^{p-1}\gamma_\phi(|\bar x-\bar y|)
      -C|\bar x-\bar y|^\sigma. \qedhere
\end{align*}
\end{proof}

\begin{lemma}\label{lem:farmodulus}
Assume that $K$ satisfies \eqref{eq:kernel} and that
$\norm{u}_{L^\infty(B_{3/2})}+\Tail(u;0,3/2)\le1$.
For every pair $(\bar x,\bar y)$ satisfying \eqref{eq:max}, the far-field integral satisfies
\begin{align}
 |I_{\mathrm{far}}|&\le C_{\mathrm{far}},\notag\\
 I_{\mathrm{far}}&\ge -C_{\mathrm{far}}\bigl(\omega_\rho(k;|\bar x-\bar y|)+|\bar x-\bar y|\bigr).
 \label{eq:farcomparison}
\end{align}
Here $C_{\mathrm{far}}$ depends only on $n,p,\Lambda$, and $\omega_\rho(k;\cdot)$
is defined in \eqref{eq:kernelmodulus}.
\end{lemma}
\begin{proof}
We have $v_{\bar x}(\bar x)=u(\bar x)$ and $v_{\bar y}(\bar y)=u(\bar y)$.
For $|z|\ge\rho>\delta$, the points $\bar x+z$ and $\bar y+z$ lie outside
$B_\delta(\bar x)$ and $B_\delta(\bar y)$, respectively, so
\[
 \begin{aligned}
 \Delta_z v_{\bar x}(\bar x)&=u(\bar x)-u(\bar x+z),\\
 \Delta_z v_{\bar y}(\bar y)&=u(\bar y)-u(\bar y+z).
 \end{aligned}
\]
The triangle inequality therefore gives
\[
 |I_{\mathrm{far}}|\le\int_{|z|\ge\rho}
 \bigl(|u(\bar x)-u(\bar x+z)|^{p-1}+|u(\bar y)-u(\bar y+z)|^{p-1}\bigr)K(z)\,\mathrm dz.
\]
In the first term, write $w=\bar x+z$. By \eqref{eq:normalized},
$|u(\bar x)|\le1$ and $|u(\bar x)-u(w)|\le2$ for $w\in B_{3/2}$.
For $|w|\ge3/2$, we have $|u(\bar x)-u(w)|^{p-1}\le C(1+|u(w)|^{p-1})$
and $|w-\bar x|\ge |w|-11/8\ge |w|/12$. Hence \eqref{eq:kernel} gives
\[
 \begin{aligned}
 \int_{|z|\ge\rho}|u(\bar x)-u(\bar x+z)|^{p-1}K(z)\,\mathrm dz
 &\le C\int_{|z|\ge\rho}\frac{\mathrm dz}{|z|^{n+sp}}
     +C\int_{|w|\ge3/2}\frac{|u(w)|^{p-1}}{|w|^{n+sp}}\,\mathrm dw\\
 &\le C\rho^{-sp}+C(3/2)^{-sp}\Tail(u;0,3/2)^{p-1}.
 \end{aligned}
\]
The same bound holds with $\bar y$ in place of $\bar x$.
Since $\rho=1/16$ and $\Tail(u;0,3/2)\le1$, we obtain $|I_{\mathrm{far}}|\le C_{\mathrm{far}}$.

For the lower bound compare the two kernels in the same integration
variable, as in \cite[Lemma~6.3]{BT}. Put
$K_x(w)=\one_{\{|w-x|\ge\rho\}}K(w-x)$ and define $K_y$ analogously.
At the contact points $u(\bar x)>u(\bar y)$ by \eqref{eq:maxproperties}.
Monotonicity of $J_p$ and positivity of $K_{\bar x}$ give
\begin{align*}
 I_{\mathrm{far}}
 &=\int_{\R^n}\bigl[J_p(u(\bar x)-u(w))K_{\bar x}(w)
                    -J_p(u(\bar y)-u(w))K_{\bar y}(w)\bigr]\,\mathrm dw\\
 &\ge\int_{\R^n}J_p(u(\bar y)-u(w))(K_{\bar x}(w)-K_{\bar y}(w))\,\mathrm dw .
\end{align*}
For every $u\in L^\infty(B_{3/2})\cap L^{p-1}_{sp}(\R^n)$ and
$x,y\in B_{11/8}$ we claim that
\begin{equation}\label{eq:farkerneldifference}
 \begin{split}
 &\int_{\R^n}
   \left(\norm{u}_{L^\infty(B_{3/2})}^{p-1}+|u(w)|^{p-1}\right)
   |K_x(w)-K_y(w)|\,\mathrm dw\\
 &\qquad\le C\bigl(\omega_\rho(k;|x-y|)+|x-y|\bigr)
   \left(\norm{u}_{L^\infty(B_{3/2})}^{p-1}+\Tail(u;0,3/2)^{p-1}\right).
 \end{split}
\end{equation}

On the common support of $K_x$ and $K_y$, both $|w-x|$ and $|w-y|$ are
at least $\rho$. Since $x,y\in B_{11/8}$, the triangle inequality gives
\[
 1+|w|\le 1+11/8+\min\{|w-x|,|w-y|\}
 \le C\min\{|w-x|,|w-y|\}.
\]
Since $|(w-x)-(w-y)|=|x-y|$, \eqref{eq:kernelmodulus} gives
$|k(w-x)-k(w-y)|\le\omega_\rho(k;|x-y|)$. The mean value theorem for
$r\mapsto r^{-n-sp}$ then gives
\[
 |K(w-x)-K(w-y)|
 \le C\left(
 \frac{\omega_\rho(k;|x-y|)}{(1+|w|)^{n+sp}}
 +\frac{|x-y|}{(1+|w|)^{n+sp+1}}\right).
\]
Since $|u(w)|\le\norm{u}_{L^\infty(B_{3/2})}$ for a.e.\ $w\in B_{3/2}$,
splitting the integral at $|w|=3/2$ and using \eqref{eq:tail} gives
\[
 \begin{aligned}
 \int_{\R^n}\frac{\norm{u}_{L^\infty(B_{3/2})}^{p-1}+|u(w)|^{p-1}}{(1+|w|)^{n+sp}}\,\mathrm dw
 &\le 2\norm{u}_{L^\infty(B_{3/2})}^{p-1}\int_{\R^n}\frac{\mathrm dw}{(1+|w|)^{n+sp}}
     +\int_{|w|\ge3/2}\frac{|u(w)|^{p-1}}{|w|^{n+sp}}\,\mathrm dw\\
 &\le C\norm{u}_{L^\infty(B_{3/2})}^{p-1}+(3/2)^{-sp}\Tail(u;0,3/2)^{p-1}.
 \end{aligned}
\]
The remaining contribution comes from $B_\rho(x)\mathbin{\triangle}B_\rho(y)
\subset B_{3/2}$, where $|K_x-K_y|\le\Lambda\rho^{-n-sp}$.
If $w\in B_\rho(x)\setminus B_\rho(y)$, then $\rho\le |w-y|\le |w-x|+|x-y|$, so
\[
 (\rho-|x-y|)_+\le |w-x|<\rho.
\]
Interchanging $x$ and $y$ gives the same shell bound for $B_\rho(y)\setminus B_\rho(x)$.
Adding the two shell volumes gives
\[
 \begin{aligned}
 |B_\rho(x)\mathbin{\triangle}B_\rho(y)|
 &\le 2n|B_1|\int_{(\rho-|x-y|)_+}^{\rho} t^{n-1}\,\mathrm dt\\
 &\le 2n|B_1|\rho^{n-1}|x-y|.
 \end{aligned}
\]
On this symmetric difference, $|u(w)|\le\norm{u}_{L^\infty(B_{3/2})}$ almost everywhere.
Hence
\[
 \begin{aligned}
 &\int_{B_\rho(x)\mathbin{\triangle}B_\rho(y)}
 \bigl(\norm{u}_{L^\infty(B_{3/2})}^{p-1}+|u(w)|^{p-1}\bigr)
 |K_x(w)-K_y(w)|\,\mathrm dw\\
 &\quad\le 2\norm{u}_{L^\infty(B_{3/2})}^{p-1}\Lambda\rho^{-n-sp}
 |B_\rho(x)\mathbin{\triangle}B_\rho(y)|\\
 &\quad\le 4n|B_1|\Lambda\rho^{-1-sp}|x-y|\norm{u}_{L^\infty(B_{3/2})}^{p-1}.
 \end{aligned}
\]
Since $\rho=1/16$, this proves \eqref{eq:farkerneldifference}.
Applying \eqref{eq:farkerneldifference} with $x=\bar x$ and $y=\bar y$ and using
\eqref{eq:normalized} gives \eqref{eq:farcomparison}.
\end{proof}

\section{Proofs of the main theorems}\label{sec:estimates}

The equation and Lemma~\ref{lem:farmodulus} give an upper bound for
$I_{B_\rho}$ at a positive maximum. A constant bound gives the bounded-data
estimate. For the refined estimate we retain the oscillations of $f$ and $k$
at each scale. If $k$ is Dini continuous outside $B_\rho$, its oscillation
contributes a bounded integral.

\subsection{Construction of the concave modulus \texorpdfstring{$\phi$}{φ}}\label{subsec:comparison-profile}

Let $p>1$ and let $\omega:(0,\rho]\to[0,\infty)$ be bounded and nondecreasing.
At a positive maximum where $I_{B_\rho}\le\omega(|\bar x-\bar y|)$,
the near-field estimate \eqref{eq:vanishing} gives
\[
 cL^{p-1}\bigl(\phi'(|\bar x-\bar y|)\bigr)^{p-1}\gamma_\phi(|\bar x-\bar y|)
 \le \omega(|\bar x-\bar y|)+C|\bar x-\bar y|^{\sigma}.
\]
To make the left-hand side dominate the right-hand side for large $L$,
we prescribe a lower bound for $(\phi')^{p-1}\gamma_\phi$.
In the variable $\tau=1+\log(\rho/r)$, the identity $\gamma_\phi(r)=-r\phi''(r)/\phi'(r)$ gives
\[
 \frac{\mathrm d}{\mathrm d\tau}\bigl(\phi'(r)\bigr)^{p-1}
 =(p-1)\bigl(\phi'(r)\bigr)^{p-1}\gamma_\phi(r).
\]
We therefore construct $(\phi'(r))^{p-1}$ by integrating a continuous
majorant of $g(\tau)=\omega(\rho e^{1-\tau})$, $\tau\ge1$. Extend $g$ by
$g(1)$ to $\tau<1$ and take its one-sided average
\begin{equation}\label{eq:qregularization}
 \bar g(\tau)=\int_0^\infty e^{-\theta}g(\tau-\theta)\,\mathrm d\theta
 =e^{-\tau}\int_{-\infty}^{\tau}e^\theta g(\theta)\,\mathrm d\theta,
 \qquad \tau\ge1.
\end{equation}
Since $g$ is nonincreasing, $g(\tau)\le\bar g(\tau)\le g(1)$.
The second expression shows that $\bar g$ is locally absolutely continuous.
We add $\tau^{-2}$ to this average because it dominates $r^\sigma$ near
$r=0$ and has finite integral over $[1,\infty)$.

Set $\phi(0)=0$ and prescribe its positive derivative by
\begin{equation}\label{eq:dataprofile}
 (\phi'(r))^{p-1}
 =1+(p-1)\int_1^\tau\bigl(\bar g(\theta)+\theta^{-2}\bigr)\,\mathrm d\theta,
 \qquad \tau=1+\log(\rho/r).
\end{equation}
Since $\bar g(\theta)+\theta^{-2}\le1+\omega(\rho)$, we have
$1\le(\phi'(r))^{p-1}\le C\tau$. Hence $\phi(r)=\int_0^r\phi'(t)\,\mathrm dt$
is finite and continuous at zero. On $(0,\rho)$, $\phi$ is $C^2$, strictly
increasing and concave, with $\phi'\ge1$. Extend $\phi$ continuously and
monotonically to $[0,3]$.

\begin{lemma}\label{lem:profileconcavity}
Let $p>1$ and let $\omega:(0,\rho]\to[0,\infty)$ be bounded and nondecreasing.
Let $\phi$ be defined by \eqref{eq:dataprofile}. With $\tau=1+\log(\rho/r)$, we have
\begin{equation}\label{eq:datacurvature}
 (\phi'(r))^{p-1}\gamma_\phi(r)=\bar g(\tau)+\tau^{-2}
 \ge\omega(r)+\tau^{-2},\qquad 0<r<\rho.
\end{equation}
Moreover,
\begin{equation}\label{eq:dataprofilebounds}
 C^{-1}\gamma_\phi(r)\le\gamma_\phi(t)\le C\gamma_\phi(r),
 \qquad r\in(0,8\rho/9),\quad t\in[7r/8,9r/8],
\end{equation}
where $C\ge1$ depends only on $p$ and $\omega(\rho)$.
\end{lemma}
\begin{proof}
Differentiating \eqref{eq:dataprofile} with respect to $\tau$ and using
$\bar g(\tau)\ge g(\tau)=\omega(r)$ gives \eqref{eq:datacurvature}.
Differentiating \eqref{eq:qregularization} gives, almost everywhere on $[1,\infty)$,
\begin{equation}\label{eq:qbounds}
 \bar g'=g-\bar g,\qquad
 -2\le \frac{\mathrm d}{\mathrm d\tau}\log\bigl(\bar g(\tau)+\tau^{-2}\bigr)\le0.
\end{equation}
Since $\phi'\ge1$ and $\bar g(\tau)\le\omega(\rho)$, \eqref{eq:datacurvature} gives
$0<\gamma_\phi(r)\le1+\omega(\rho)$. Using
$\mathrm d\log\phi'(r)/\mathrm d\tau=\gamma_\phi(r)$ in
\eqref{eq:datacurvature}, we obtain from \eqref{eq:qbounds} that
$\log\gamma_\phi(r)$ is Lipschitz in $\tau$. Thus, for $0<r<8\rho/9$
and $7/8\le t/r\le9/8$,
\[
 \left|\log\frac{\gamma_\phi(t)}{\gamma_\phi(r)}\right|
 \le C\left|\log\frac{\rho}{t}-\log\frac{\rho}{r}\right|
 = C\left|\log\frac{t}{r}\right|\le C\log\frac87.
\]
Exponentiating gives \eqref{eq:dataprofilebounds}.
\end{proof}

\begin{lemma}\label{lem:dataprofile}
Let $p>1$ and let $\omega:(0,\rho]\to[0,\infty)$ be bounded and nondecreasing.
Let $\phi$ be defined by \eqref{eq:dataprofile}. Then
\begin{equation}\label{eq:profilemodulus}
 \phi(r)\asymp r\left(1+\int_r^\rho\omega(t)\frac{\mathrm dt}{t}\right)^{1/(p-1)},
 \qquad 0<r\le\rho,
\end{equation}
with comparison constants depending only on $p$ and $\omega(\rho)$.
\end{lemma}
\begin{proof}
Integrating the first identity in \eqref{eq:qbounds} gives
\[
 0\le\int_1^\tau \bigl(\bar g(\theta)-g(\theta)\bigr)\,\mathrm d\theta
   =g(1)-\bar g(\tau)\le g(1).
\]
Since $\int_1^\tau \theta^{-2}\,\mathrm d\theta\le1$, the substitution
$t=\rho e^{1-\theta}$ in the integral of $g$ yields
\[
 0\le (\phi'(r))^{p-1}-1-(p-1)\int_r^\rho\omega(t)\frac{\mathrm dt}{t}
 \le(p-1)(1+\omega(\rho)).
\]
By \eqref{eq:dataprofile}, increasing $\tau$ by $\theta\ge0$ increases
$(\phi'(r))^{p-1}$ by at most $(p-1)(1+\omega(\rho))\theta$. Since
$\phi'(r)\ge1$, substituting $t=re^{-\theta}$ gives
\begin{align*}
 r\phi'(r)\le\phi(r)
 &=r\int_0^\infty e^{-\theta}\phi'(re^{-\theta})\,\mathrm d\theta\\
 &\le r\phi'(r)
        \int_0^\infty e^{-\theta}(1+C\theta)^{1/(p-1)}\,\mathrm d\theta
  \le Cr\phi'(r).
\end{align*}
Combining these estimates proves \eqref{eq:profilemodulus}.
\end{proof}

\subsection{Interior estimates}\label{subsec:kernel-estimates}

Given a bounded nondecreasing $\omega:(0,\rho]\to[0,\infty)$, take $\phi$
in \eqref{eq:max} to be the profile constructed from $\omega$ in
Section~\ref{subsec:comparison-profile}. We seek $L$ such that
\[
 \max_{x,y\in\overline B_{11/8}}\{u(x)-u(y)-L\phi(|x-y|)-\psi(x)\}\le0.
\]
With $c,C,\sigma$ from \eqref{eq:vanishing}, fix a sufficiently large
$L\ge4/\rho$ such that
\[
 cL^{p-1}>1+Ct^\sigma\bigl(1+\log(\rho/t)\bigr)^2,\qquad 0<t<\rho/2.
\]
The right-hand side is bounded, so $L$ can be chosen in terms of
$n,p,\lambda,\Lambda$ and $\omega(\rho)$ alone.
For $u$ satisfying \eqref{eq:normalized}, every positive maximum satisfies
$L\phi(|\bar x-\bar y|)<2\le L\phi(\rho/2)$ by \eqref{eq:maxproperties}
and $\phi(\rho/2)\ge\rho/2$. Thus $|\bar x-\bar y|<\rho/2$.
We now apply the near-field estimate of Lemma~\ref{lem:nearfield} to the
profile constructed in Section~\ref{subsec:comparison-profile}.

\begin{lemma}\label{lem:profilegain}
Assume that $u\in C(B_{3/2})\cap L^{p-1}_{sp}(\R^n)$ satisfies \eqref{eq:normalized}.
Let $\omega:(0,\rho]\to[0,\infty)$ be bounded and nondecreasing.
Let $\phi$ be defined by \eqref{eq:dataprofile}.
At every positive maximum $(\bar x,\bar y)$ in \eqref{eq:max}, the near-field integral satisfies
\begin{equation}\label{eq:profilegain}
 I_{B_\rho}>\omega(|\bar x-\bar y|).
\end{equation}
\end{lemma}
\begin{proof}
Since $\phi'\ge1$ and $\bar g\le\omega(\rho)$, \eqref{eq:datacurvature} gives
$0<\gamma_\phi\le1+\omega(\rho)$. Together with \eqref{eq:dataprofilebounds},
this verifies \eqref{eq:profileconditions} with $M_\phi\le C(p,\omega(\rho))$.
At $t=|\bar x-\bar y|<\rho/2$, the choice of $L$ gives
$cL^{p-1}>1+Ct^\sigma\tau^2$, where $\tau=1+\log(\rho/t)$.
Thus \eqref{eq:vanishing} and \eqref{eq:datacurvature} give
\[
 I_{B_\rho}-\omega(|\bar x-\bar y|)
 \ge (cL^{p-1}-1)\omega(|\bar x-\bar y|)
       +cL^{p-1}\tau^{-2}-Ct^\sigma>\tau^{-2}>0. \qedhere
\]
\end{proof}

\begin{lemma}\label{lem:normalizedmodulus}
Assume that $u\in C(B_{3/2})\cap L^{p-1}_{sp}(\R^n)$ satisfies \eqref{eq:normalized}.
Let $\omega:(0,\rho]\to[0,\infty)$ be bounded and nondecreasing.
Let $\phi$ be defined by \eqref{eq:dataprofile}.
For the constant $L\ge4/\rho$ fixed at the beginning of Section~\ref{subsec:kernel-estimates},
suppose that
\[
 \max_{x,y\in\overline B_{11/8}}\{u(x)-u(y)-L\phi(|x-y|)-\psi(x)\}\le0.
\]
Then
\begin{equation}\label{eq:normalizedmodulus}
 |u(x)-u(y)|\le C|x-y|
 \left(1+\int_{|x-y|}^\rho\omega(t)\frac{\mathrm dt}{t}\right)^{1/(p-1)},
 \qquad x,y\in B_1,\quad 0<|x-y|\le\rho,
\end{equation}
where $C$ depends only on $n,p,\lambda,\Lambda$ and $\omega(\rho)$.
\end{lemma}
\begin{proof}
Since $\psi=0$ on $B_1$, we have $u(x)-u(y)\le L\phi(|x-y|)$ for all
$x,y\in B_1$. Interchanging $x$ and $y$ gives
$|u(x)-u(y)|\le L\phi(|x-y|)$, and Lemma~\ref{lem:dataprofile} yields
\eqref{eq:normalizedmodulus}.
\end{proof}

In Corollary~\ref{cor:general-log} and Theorem~\ref{thm:kernel-modulus},
we choose $\omega$ so that the equation and the far-field estimate give
$I_{B_\rho}\le\omega(|\bar x-\bar y|)$ at any positive maximum.
Lemma~\ref{lem:profilegain} rules out such a maximum, so
Lemma~\ref{lem:normalizedmodulus} applies.

\begin{corollary}\label{cor:general-log}
Suppose that
$u\in C(B_2)\cap L^{p-1}_{sp}(\R^n)$ satisfies $-F\le\LL u\le F$ in
$B_{3/2}$ in the viscosity sense for some $F\ge0$. Then, for all
distinct $x,y\in B_1$,
\begin{equation}\label{eq:viscosity-log}
 \begin{split}
 |u(x)-u(y)|\le{}&C\left(\norm{u}_{L^\infty(B_{3/2})}+\Tail(u;0,3/2)
                   +F^{1/(p-1)}\right)\\
 &\qquad\times |x-y|\bigl(1+\bigl|\log|x-y|\bigr|\bigr)^{1/(p-1)},
 \end{split}
\end{equation}
where $C$ depends only on $n,p,\lambda,\Lambda$.
\end{corollary}
\begin{proof}
Normalize $u$ as in \eqref{eq:normalized}. By \eqref{eq:differenceineq}
and Lemma~\ref{lem:farmodulus}, $I_{B_\rho}\le 2-I_{\mathrm{far}}\le 2+C_{\mathrm{far}}$
at every positive maximum. Thus \eqref{eq:normalizedmodulus} with
$\omega(t)=2+C_{\mathrm{far}}$ gives
\[
 |u(x)-u(y)|\le C|x-y|\left(1+\log\frac{\rho}{|x-y|}\right)^{1/(p-1)}
 \qquad(0<|x-y|\le\rho,\ x,y\in B_1).
\]
For $|x-y|>\rho$, we have $|u(x)-u(y)|\le2\le32|x-y|$.
Multiplying by $A$ gives \eqref{eq:viscosity-log}.
\end{proof}

By Lemma~\ref{lem:weakregularity} and Proposition~\ref{prop:weakviscosity},
Corollary~\ref{cor:general-log} applies to every local weak solution of
$\LL u=f$ in $B_2$ with $f\in L^\infty_{\mathrm{loc}}(B_2)$ and
$F=\norm{f}_{L^\infty(B_{3/2})}$. For the standard kernel, $\lambda=\Lambda=1$,
which gives Theorem~\ref{thm:log}.

\Needspace{4\baselineskip}
\begin{theorem}\label{thm:kernel-modulus}
Assume that
$\int_0^\rho\omega_\rho(k;t)\,\mathrm dt/t<\infty$.
Let $f\in L^\infty_{\mathrm{loc}}(B_2)$ and
$u\in C(B_2)\cap L^{p-1}_{sp}(\R^n)$. Suppose that on every ball
$B\Subset B_2$ the inequalities
\[
 \operatorname*{ess\,inf}_B f\le\LL u\le\operatorname*{ess\,sup}_B f
\]
hold in the viscosity sense. Then, for all $x,y\in B_1$ with $0<|x-y|\le1$,
\begin{equation}\label{eq:datumestimate}
 \begin{split}
 |u(x)-u(y)|\le C|x-y|\Bigg[&\norm{u}_{L^\infty(B_{3/2})}+\Tail(u;0,3/2)\\
 &+\left(\norm{f}_{L^\infty(B_{3/2})}
       +\int_{|x-y|}^1\omega(f;t)\frac{\mathrm dt}{t}
   \right)^{1/(p-1)}\Bigg],
 \end{split}
\end{equation}
where $\omega(f;t)$ is defined in \eqref{eq:omega} and $C$ depends only on
$n,p,\lambda,\Lambda$ and
$\int_0^\rho\omega_\rho(k;t)\,\mathrm dt/t$.
\end{theorem}
\begin{proof}
Taking $B=B_{3/2}$, we normalize $u$ as in \eqref{eq:normalized} with
\[
 A=\norm{u}_{L^\infty(B_{3/2})}+\Tail(u;0,3/2)
   +\norm{f}_{L^\infty(B_{3/2})}^{1/(p-1)}.
\]
The normalized datum $\widehat f=A^{1-p}f$ satisfies $\norm{\widehat f}_{L^\infty(B_{3/2})}\le1$.
Take
\[
 \omega(t)=\omega(\widehat f;2t)+C_{\mathrm{far}}\bigl(\omega_\rho(k;t)+t\bigr),
 \qquad 0<t\le\rho.
\]
Since $\omega(\widehat f;2t)\le2$ and
$\omega_\rho(k;t)\le\Lambda-\lambda$, $\omega(\rho)$ is bounded in terms of $n,p,\Lambda$.

At a positive maximum $(\bar x,\bar y)$ in \eqref{eq:max}, take
$B=B_{|\bar x-\bar y|}((\bar x+\bar y)/2)\Subset B_{3/2}$.
Since $\delta<|\bar x-\bar y|/8$, both test balls are compactly contained in $B$.
Applying the viscosity inequalities on $B$ at the two contact points gives
\[
 \LL v_{\bar x}(\bar x)-\LL v_{\bar y}(\bar y)
 \le \operatorname*{ess\,sup}_B\widehat f-\operatorname*{ess\,inf}_B\widehat f
 \le \omega(\widehat f;2|\bar x-\bar y|).
\]
The operator difference equals $I_{B_\rho}+I_{\mathrm{far}}$ by
\eqref{eq:differenceineq}, so \eqref{eq:farcomparison} gives
$I_{B_\rho}\le\omega(|\bar x-\bar y|)$, contradicting \eqref{eq:profilegain}.
For $0<r\le\rho$, the two contributions in \eqref{eq:normalizedmodulus} satisfy
\begin{align*}
 \int_r^\rho\omega(\widehat f;2t)\frac{\mathrm dt}{t}
 &=\int_{2r}^{2\rho}\omega(\widehat f;t)\frac{\mathrm dt}{t}
 \le\int_r^1\omega(\widehat f;t)\frac{\mathrm dt}{t},\\
 \int_r^\rho\bigl(\omega_\rho(k;t)+t\bigr)\frac{\mathrm dt}{t}
 &\le\int_0^\rho\omega_\rho(k;t)\frac{\mathrm dt}{t}+\rho.
\end{align*}
For $\rho<|x-y|\le1$, we use $|u(x)-u(y)|\le2\le32|x-y|$.
Multiplying by $A$ gives \eqref{eq:datumestimate}, since
$\omega(\widehat f;t)=A^{1-p}\omega(f;t)$.
\end{proof}

By Lemma~\ref{lem:weakregularity} and Proposition~\ref{prop:weakviscosity},
\eqref{eq:datumestimate} applies to every local weak solution of $\LL u=f$
in $B_2$ with $f\in L^\infty_{\mathrm{loc}}(B_2)$. For the standard kernel,
$\lambda=\Lambda=1$ and $\omega_\rho(k;t)=0$ for all $t>0$, which gives
Theorem~\ref{thm:logmodulus} with $C=C(n,p)$.

\begin{remark}\label{rem:kerneloscillation}
Assume that $u$ and $f$ satisfy the hypotheses
on the solution and the datum in Theorem~\ref{thm:kernel-modulus}.
The following estimate holds even when $\int_0^\rho\omega_\rho(k;t)\,\mathrm dt/t=\infty$.
Set
$M=\norm{u}_{L^\infty(B_{3/2})}+\Tail(u;0,3/2)$.
For $x,y\in B_1$ with $0<|x-y|\le1$, we have
\[
 \begin{aligned}
 \frac{|u(x)-u(y)|}{|x-y|}
 &\le CM\left(1+\int_{|x-y|}^1\omega_\rho(k;t)\frac{\mathrm dt}{t}
        \right)^{1/(p-1)}\\
 &\quad+C\left(\norm{f}_{L^\infty(B_{3/2})}
        +\int_{|x-y|}^1\omega(f;t)\frac{\mathrm dt}{t}\right)^{1/(p-1)},
 \end{aligned}
\]
where $C=C(n,p,\lambda,\Lambda)$.
Indeed, for $u/A$ the right-hand side of \eqref{eq:farkerneldifference}
carries the factor $(M/A)^{p-1}$ up to a constant depending only on $p$.
Keeping this factor in the kernel term of $\omega$ in \eqref{eq:normalizedmodulus} and multiplying
by $A=M+\norm{f}_{L^\infty(B_{3/2})}^{1/(p-1)}$ gives the displayed estimate.
\end{remark}

In Theorem~\ref{thm:kernel-modulus}, a Dini modulus in
\eqref{eq:datumestimate} yields a Lipschitz estimate, and a logarithmic
modulus yields the upper bounds in Remark~\ref{rem:sharp-endpoints}.

\section{Sharpness for a prescribed data modulus}\label{sec:sharpness}

For each prescribed modulus, we construct a solution that attains the
bound in Theorem~\ref{thm:logmodulus} up to constants at every sufficiently
small scale. Here $\omega(f;r)$ is the modulus of $f$ on $B_{3/2}$
defined in \eqref{eq:omega}.

\begin{theorem}\label{thm:sharp-modulus}
Let $n\ge1$, $1<p<\infty$ and $s=(p-1)/p$. For each nondecreasing,
subadditive $\omega:(0,1]\to[0,\infty)$ with $\omega(1)=1$, there exist
$u\in C(\R^n)\cap L^\infty(\R^n)\cap W^{s,p}_{\mathrm{loc}}(\R^n)$
and $f\in L^\infty(B_2)$
such that $(-\Delta_p)^s u=f$ weakly in $B_2$, $u(0)=0$, and
$\norm{u}_{L^\infty(\R^n)}+\norm{f}_{L^\infty(B_2)}\le C$.
For every $0<r<r_*$, where $r_*\in(0,1)$ depends only on $n,p$,
the same solution satisfies
\begin{align}
 \omega(f;r)&\asymp\omega(r),\label{eq:sharp-data}\\
 u(re_1)-u(0)&\asymp
 r\left(1+\int_r^1\omega(t)\frac{\mathrm dt}{t}\right)^{1/(p-1)},
 \label{eq:sharp-increment}
\end{align}
where $e_1$ is the first coordinate vector. The constants in these
estimates depend only on $n,p$. If $\omega(0+)=0$, then $f$ is
continuous in $B_2$.
\end{theorem}

Theorem~\ref{thm:sharp-modulus} applies to every nonzero data modulus
in the sense of \eqref{eq:omega}. For continuous data this is the usual
modulus of continuity. The construction also allows $\omega(0+)>0$.
To normalize a prescribed modulus, set $M=\omega(1)>0$ and apply the
theorem to $\omega/M$. If $(u,f)$ is the resulting pair, then
$M^{1/(p-1)}u$ has right-hand side $Mf$, whose modulus is comparable
to $\omega$ at every sufficiently small scale.

The modulus in \eqref{eq:omega} is nondecreasing since the set of
admissible pairs grows with the radius. To see subadditivity, fix
$r,t>0$ with $r+t\le1$. For $x,y\in B_{3/2}$ with $|x-y|\le r+t$,
the point $z=(tx+ry)/(r+t)$ lies in $B_{3/2}$ by convexity and
satisfies $|x-z|\le r$ and $|z-y|\le t$. The maps
$(x,y)\mapsto(x,z)$ and $(x,y)\mapsto(z,y)$ are invertible and linear,
so the exceptional sets in \eqref{eq:omega} have null preimages.
Thus for almost every such pair $(x,y)$,
\[
 |f(x)-f(y)|\le |f(x)-f(z)|+|f(z)-f(y)|
 \le \omega(f;r)+\omega(f;t).
\]
Taking the essential supremum gives $\omega(f;r+t)\le\omega(f;r)+\omega(f;t)$.

\subsection{Construction of the one-dimensional profile}\label{subsec:sharp-profile}

We prescribe the slope of an odd function $v$ using
$q(t)=\Omega(t)-t\Omega'(t)$, where $\Omega$ is smooth, concave and
comparable to $\omega$. Identity \eqref{eq:sharp-primitive} determines
the growth of the slope, while \eqref{eq:sharp-hardy} controls the
contribution of displacements $z>2r$ to the operator.

\begin{lemma}\label{lem:sharp-concave}
For $\omega$ as in Theorem~\ref{thm:sharp-modulus}, there is a smooth,
nondecreasing, concave function $\Omega$ on $(0,\infty)$ which is
comparable to $\omega$ on $(0,1]$, equals $1$ near $[1,\infty)$ and
satisfies $\Omega(0+)=\omega(0+)$. The function
$q(r)=\Omega(r)-r\Omega'(r)$ satisfies
\begin{equation}\label{eq:sharp-qbounds}
 0\le q(r)\le\Omega(r)\le1,\qquad 0\le rq'(r)\le C\Omega(r)
 \qquad(0<r\le1).
\end{equation}
Moreover,
\begin{align}
 \int_r^1 q(t)\,\frac{\mathrm dt}{t}
 &=\int_r^1\Omega(t)\,\frac{\mathrm dt}{t}-1+\Omega(r),
      &&0<r\le1,\label{eq:sharp-primitive}\\
 \Omega(r)&=\frac rR\Omega(R)+r\int_r^R\frac{q(t)}{t^2}\,\mathrm dt,
      &&0<r<R\le1.\label{eq:sharp-hardy}
\end{align}
All comparison constants here are absolute.
\end{lemma}
\begin{proof}
Set $\omega(0)=0$. Its least concave majorant $\Omega_0$ on $[0,1]$
is the supremum of $\sum_j\lambda_j\omega(t_j)$ over finite convex
combinations with $\sum_j\lambda_jt_j=r$. For $0<r\le1$ and
$0\le t\le1$, monotonicity and subadditivity give
$\omega(t)\le(1+t/r)\omega(r)$. Averaging in each convex combination
yields $\omega(r)\le\Omega_0(r)\le2\omega(r)$. Moving the points $t_j$
towards $1$ increases their mean without decreasing the sum, so
$\Omega_0$ is nondecreasing. Its endpoint values are $0$ and $1$.
For $0<r<d\le1$, the total weight assigned to $t_j>d$ is at most
$r/d$, hence
\[
 \omega(r)\le\Omega_0(r)\le\omega(d)+r/d.
\]
Taking first $r\downarrow0$ and then $d\downarrow0$ gives the asserted
right limit for $\Omega_0$.

Extend $\Omega_0$ by $1$ on $[1,\infty)$. Fix
$0\le\eta\in C_c^\infty((1,2))$ with $\int_1^2\eta(\sigma)\,\mathrm d\sigma=1$ and set
\[
 \Omega(r)=\int_1^2\eta(\sigma)\Omega_0(\sigma r)\,\mathrm d\sigma.
\]
This average preserves monotonicity and concavity. Since
$\Omega_0(r)\le\Omega(r)\le\Omega_0(2r)\le2\Omega_0(r)$, it preserves
comparability and the right limit at zero. In the variable $\log r$
it is a convolution with a fixed smooth kernel. Differentiating this convolution kernel once and twice gives
\[
 r|\Omega'(r)|+r^2|\Omega''(r)|\le C\Omega_0(2r)\le C\Omega(r).
\]
Concavity gives $0\le r\Omega'(r)\le\Omega(r)$ and
$q'(r)=-r\Omega''(r)\ge0$. This proves \eqref{eq:sharp-qbounds}.
Integrating $q(t)/t=\Omega(t)/t-\Omega'(t)$ and
$(\Omega(t)/t)'=-q(t)/t^2$ proves the two identities.
\end{proof}

Fix a smooth nonincreasing function $\chi$ that is
$1$ on $[0,2]$, positive on $[0,4)$ and zero on $[4,\infty)$. Define
\begin{equation}\label{eq:sharp-slope}
 a(r)=\chi(r)\left(8+\int_r^1 q(t)\,\frac{\mathrm dt}{t}\right)^{1/(p-1)}
 \quad(0<r<4),\qquad a(r)=0\quad(r\ge4).
\end{equation}
Set $v(r)=\int_0^r a(t)\,\mathrm dt$ for $r\ge0$ and extend $v$ oddly
to $\R$. Since $q=1$ on $[1,\infty)$, the expression inside the power in
\eqref{eq:sharp-slope} is at least $8-\log4$. The slope $a$ is nonincreasing.
Thus $v$ is
nondecreasing on $\R$, concave on $[0,\infty)$ and constant on each of
the half-lines $|x|\ge4$. In the one-dimensional estimates below,
$c,C>0$ depend only on $p$, with $\eta$ and $\chi$ fixed.

\begin{lemma}\label{lem:sharp-profile}
For $0<r\le1$ the slope and its primitive satisfy
\begin{equation}\label{eq:sharp-growth}
 a(r)^{p-1}\asymp1+\int_r^1\Omega(t)\,\frac{\mathrm dt}{t},\qquad
 ra(r)\le v(r)\le C ra(r),\qquad a(r)\le C\frac{\Omega(r)}r.
\end{equation}
The variation of the slope satisfies
\begin{align}
 -a'(r)&=\frac{q(r)a(r)^{2-p}}{(p-1)r},&
 |a''(r)|&\le C\frac{\Omega(r)a(r)^{2-p}}{r^2},\label{eq:sharp-derivatives}\\
 \int_0^r(a(t)-a(r))\,\mathrm dt
 &\le C r q(r)a(r)^{2-p}.&&\label{eq:sharp-average}
\end{align}
Moreover, $\norm{v}_{L^\infty(\R)}\le C$ and
$v\in W^{1,p}_{\mathrm{loc}}(\R)$.
\end{lemma}
\begin{proof}
By \eqref{eq:sharp-primitive},
$a(r)^{p-1}=7+\Omega(r)+\int_r^1\Omega(t)\,\mathrm dt/t$.
Using $a(r)^{p-1}\ge8$ and $a(t)^{p-1}-a(r)^{p-1}\le\log(r/t)$, we obtain
\begin{equation}\label{eq:sharp-sloperatio}
 1\le\frac{a(t)}{a(r)}
 \le\left(1+\tfrac18\log(r/t)\right)^{1/(p-1)},\qquad 0<t\le r\le1.
\end{equation}
Integration in $t$ gives $ra(r)\le v(r)\le Cra(r)$.
For $r\le t\le1$, concavity gives
$\Omega(t)\le\min\{1,t\Omega(r)/r\}$ and $\Omega(r)\ge r$.
Splitting at $t=r/\Omega(r)$ gives
$\int_r^1\Omega(t)\,\mathrm dt/t\le1+\log(\Omega(r)/r)$.
Since $(1+\log b)^{1/(p-1)}\le Cb$ for $b\ge1$, this proves
\eqref{eq:sharp-growth}.

Differentiating \eqref{eq:sharp-slope} gives
$(p-1)a^{p-2}a'=-q/r$. Differentiating once more and using
$|rq'|+q+q^2/a^{p-1}\le C\Omega$ gives \eqref{eq:sharp-derivatives}.
By Tonelli's theorem,
\[
 \int_0^r(a(t)-a(r))\,\mathrm dt
 =\frac{1}{p-1}\int_0^r q(t)a(t)^{2-p}\,\mathrm dt.
\]
Here $q(t)\le q(r)$. For $p\ge2$ we have $a(t)^{2-p}\le a(r)^{2-p}$.
For $1<p<2$, integrate \eqref{eq:sharp-sloperatio} using
$\int_0^1(1+|\log t|)^{(2-p)/(p-1)}\,\mathrm dt<\infty$.
This gives \eqref{eq:sharp-average}. Finally,
$a(r)\le(8+\log(1/r))^{1/(p-1)}$ on $(0,1)$, and $a$ is fixed on
$[1,\infty)$. Thus $a\in L^1(0,\infty)\cap L^p(0,\infty)$ with
uniform bounds, proving the last assertions.
\end{proof}

\subsection{Estimates for the one-dimensional operator}\label{subsec:sharp-data}

Let $g_\varepsilon$ denote the one-dimensional integral in \eqref{eq:fracplap}
for $v$, truncated to $|x-y|>\varepsilon$. For $0<r<2$ and $z>0$, write
\[
 d_-=v(r)-v(r-z),\qquad d_+=v(r+z)-v(r),\qquad
 N(r,z)=d_-^{p-1}-d_+^{p-1}.
\]
Both increments are positive. Since $v'(x)=a(|x|)$ almost everywhere
and $a$ is nonincreasing,
\[
 d_--d_+=\int_0^z\bigl(a(|r-t|)-a(r+t)\bigr)\,\mathrm dt\ge0.
\]
Since $1+sp=p$, pairing the displacements $z$ and $-z$ gives
\begin{equation}\label{eq:sharp-paired}
 g_\varepsilon(r)=\int_\varepsilon^\infty N(r,z)\frac{\mathrm dz}{z^p},
 \qquad N(r,z)\ge0.
\end{equation}
We will show that the principal value $g=(-\Delta_p)^s v$ satisfies
\[
 g(r)\asymp\Omega(r),\quad 0<r<1/16,\qquad
 |g'(r)|\le C\frac{\Omega(r)}r,\quad 0<r\le3/2.
\]
The first estimate controls the oscillation near zero and, by oddness,
gives $g(r/4)-g(-r/4)=2g(r/4)\ge c\Omega(r)$ at small scales.
The derivative estimate controls differences between nearby points on
the same side of zero whose distance from zero exceeds their separation.

\begin{lemma}\label{lem:sharp-small-displacements}
For $0<r<1/16$,
\[
 \int_0^{2r} N(r,z)\frac{\mathrm dz}{z^p}\le C\Omega(r),
 \qquad
 \int_0^{2r}|\partial_r N(r,z)|\frac{\mathrm dz}{z^p}
 \le C\frac{\Omega(r)}r.
\]
\end{lemma}
\begin{proof}
For $0<z<r/2$, \eqref{eq:sharp-sloperatio} gives
$d_\pm\asymp za(r)$. Integrating $a'$ between $r-z$ and $r+z$ and
using \eqref{eq:sharp-derivatives} gives
$0\le d_--d_+\le Cz^2\Omega(r)a(r)^{2-p}/r$.
The mean value theorem then gives $N\le Cz^p\Omega(r)/r$.
For the derivative we use the cancellation
\[
 \frac{1}{p-1}\partial_r N
 =(d_+^{p-2}-d_-^{p-2})(a(r)-a(r+z))
  +d_-^{p-2}\bigl(2a(r)-a(r-z)-a(r+z)\bigr).
\]
Applying \eqref{eq:sharp-derivatives} to the two terms gives
\[
 |\partial_r N(r,z)|
 \le C\frac{z^p}{r^2}
       \left(\frac{\Omega(r)^2}{a(r)^{p-1}}+\Omega(r)\right)
 \le C\frac{z^p\Omega(r)}{r^2},
\]
since $a(r)^{p-1}\ge8$ and $\Omega(r)\le1$. Integrating the bounds for
$N$ and $|\partial_r N|$ against $\mathrm dz/z^p$ over $(0,r/2)$ gives
$C\Omega(r)$ and $C\Omega(r)/r$, respectively.

For $r/2\le z\le2r$, we have $za(3r)\le d_+\le d_-\le2v(3r)$, so
$d_\pm\asymp ra(r)$. Subtracting the linear function $x\mapsto a(3r)x$
from $v$ and using \eqref{eq:sharp-average} gives
\[
 d_--d_+\le 4\int_0^{3r}\bigl(a(t)-a(3r)\bigr)\,\mathrm dt
 \le C r\Omega(r)a(r)^{2-p}.
\]
For almost every $z$ in this range,
\[
 |\partial_r N|\le C(ra(r))^{p-2}
 \bigl(|a(r)-a(|r-z|)|+|a(r)-a(r+z)|\bigr).
\]
The integral of the parentheses in $z$ is bounded by
$C\int_0^{3r}(a(t)-a(3r))\,\mathrm dt$. Indeed, $|r-z|\le r$,
and $r(a(r)-a(3r))\le\int_0^r(a(t)-a(3r))\,\mathrm dt$.
Since $z\asymp r$, the mean value theorem for $d_\pm^{p-1}$ and the
bound for $|\partial_r N|$ give
\[
 \int_{r/2}^{2r}N(r,z)\,\frac{\mathrm dz}{z^p}\le C\Omega(r),\qquad
 \int_{r/2}^{2r}|\partial_r N(r,z)|\,\frac{\mathrm dz}{z^p}
 \le C\frac{\Omega(r)}r. \qedhere
\]
\end{proof}

\begin{lemma}\label{lem:sharp-large-displacements}
For $0<r<1/16$,
\[
 \int_{2r}^\infty N(r,z)\frac{\mathrm dz}{z^p}\asymp\Omega(r),
 \qquad
 \int_{2r}^\infty|\partial_r N(r,z)|\frac{\mathrm dz}{z^p}
 \le C\frac{\Omega(r)}r.
\]
\end{lemma}
\begin{proof}
For the lower bound, set $R=1/16$ and recall from \eqref{eq:sharp-hardy} that
\[
 \Omega(r)=r\int_r^R\frac{q(t)}{t^2}\,\mathrm dt+\frac rR\Omega(R).
\]
We bound the first term using $2r<z<3R$ and the second using $5<z<6$.
For $2r<z<1/2$, the bounds $za(z+r)\le d_+\le d_-\le2v(z)$ give
$d_\pm\asymp za(z)$. Oddness gives
\[
 d_--d_+=2v(r)-\int_{z-r}^{z+r}a(t)\,\mathrm dt.
\]
By monotonicity of $a$,
\[
 d_--d_+\ge2r\bigl(a(r)-a(z-r)\bigr)
 =\frac{2r}{p-1}\int_r^{z-r}\frac{q(t)}{t a(t)^{p-2}}\,\mathrm dt.
\]
The mean value theorem gives $N\ge c(za(z))^{p-2}(d_--d_+)$.
Restricting the resulting double integral to $r\le t\le R$ and
$2t\le z\le3t$, where $z-r\ge t$ and $a(z)\asymp a(t)$, yields
\[
 \int_{2r}^{3R}N\,\frac{\mathrm dz}{z^p}
 \ge cr\int_r^R\frac{q(t)}t
   \int_{2t}^{3t}\left(\frac{a(z)}{a(t)}\right)^{p-2}
                    \frac{\mathrm dz}{z^2}\,\mathrm dt
 \ge cr\int_r^R\frac{q(t)}{t^2}\,\mathrm dt.
\]
On $5<z<6$, the function $v$ takes its constant values at $r-z$ and $r+z$,
so $d_-=v(4)+v(r)$ and $d_+=v(4)-v(r)$. Both increments lie in a fixed
compact subinterval of $(0,\infty)$. The mean value theorem gives
$N\ge c v(r)\ge cr$. Since $\Omega(R)\le1$ and $R=1/16$,
\[
 \int_5^6 N\,\frac{\mathrm dz}{z^p}\ge cr\ge\frac{c}{16}\frac rR\Omega(R).
\]
The two intervals are disjoint and $N\ge0$, so adding the bounds gives
$\int_{2r}^\infty N(r,z)\,\mathrm dz/z^p\ge c\Omega(r)$.

For the upper bounds on $2r<z<1/2$, monotonicity of $a$ gives
\[
 d_--d_+\le 2\int_0^r(a(t)-a(r))\,\mathrm dt+2r(a(r)-a(z+r)).
\]
The mean value theorem gives
$N/z^p\le Ca(z)^{p-2}(d_--d_+)/z^2$. To integrate this weight,
we use
\begin{equation}\label{eq:sharp-weight}
 \int_{t/2}^1\left(\frac{a(z)}{a(t)}\right)^{p-2}
       \frac{\mathrm dz}{z^2}\le\frac Ct,\qquad 0<t\le1.
\end{equation}
On $[t/2,t]$ the slopes are comparable. For $z\ge t$, the ratio raised
to $p-2$ is at most $1$ if $p\ge2$. If $1<p<2$,
\eqref{eq:sharp-sloperatio} bounds it by
$C(1+\log(z/t))^{(2-p)/(p-1)}$. Substituting $z=t\sigma$ proves
\eqref{eq:sharp-weight}, since every fixed power of $1+\log\sigma$ is
integrable against $\sigma^{-2}\,\mathrm d\sigma$ on $[1,\infty)$.

The term containing $\int_0^r(a(t)-a(r))\,\mathrm dt$ contributes at
most $Cq(r)$ by \eqref{eq:sharp-average} and \eqref{eq:sharp-weight}
with $t=r$. For the second term, expand $a(r)-a(z+r)$ using
\eqref{eq:sharp-derivatives}. Since $t\le z+r\le3z/2$ in this
integral, Tonelli's theorem and \eqref{eq:sharp-weight} give
\begin{align}
 \int_{2r}^{1/2} a(z)^{p-2}\bigl(a(r)-a(z+r)\bigr)\frac{\mathrm dz}{z^2}
 &\le \frac{1}{p-1}\int_r^{3/4}\frac{q(t)}t
       \left[\int_{t/2}^1\left(\frac{a(z)}{a(t)}\right)^{p-2}
             \frac{\mathrm dz}{z^2}\right]\,\mathrm dt\notag\\
 &\le C\int_r^{3/4}\frac{q(t)}{t^2}\,\mathrm dt.
 \label{eq:sharp-shift}
\end{align}
For $2r<z<1/2$,
\[
 \frac{1}{p-1}\partial_r N
 =d_-^{p-2}(a(r)-a(z-r))+d_+^{p-2}(a(r)-a(z+r)).
\]
Both terms are nonnegative, and their sum is bounded by
$C(za(z))^{p-2}(a(r)-a(z+r))$. By \eqref{eq:sharp-shift} and
\eqref{eq:sharp-hardy},
\[
 \begin{aligned}
  \int_{2r}^{1/2}N(r,z)\,\frac{\mathrm dz}{z^p}
  &\le Cq(r)+Cr\int_r^{3/4}\frac{q(t)}{t^2}\,\mathrm dt\le C\Omega(r),\\
  \int_{2r}^{1/2}|\partial_r N(r,z)|\,\frac{\mathrm dz}{z^p}
  &\le C\frac{\Omega(r)}r.
 \end{aligned}
\]

For $z\ge1/2$, both increments lie in a fixed compact subinterval
of $(0,\infty)$. Indeed $d_+\ge\int_r^{r+1/2}a(t)\,\mathrm dt\ge c$
and $d_-\le2\norm{v}_{L^\infty(\R)}$. Thus
$N\le C(d_--d_+)\le 2Cv(r)$ and $|\partial_r N|\le Ca(r)$.
Integrating over $[1/2,\infty)$ and applying \eqref{eq:sharp-growth}
gives the same two upper bounds.
\end{proof}

\begin{proposition}\label{prop:sharp-rhs}
The principal value $g=(-\Delta_p)^s v$ exists on $(-2,2)$, is odd,
and satisfies $\norm{g}_{L^\infty(-2,2)}\le C$. It is continuously
differentiable away from zero and satisfies
\begin{equation}\label{eq:sharp-rhs}
 \begin{aligned}
 c\Omega(r)\le g(r)&\le C\Omega(r), & 0<r<1/16,\\
 |g'(r)|&\le C\,\frac{\Omega(r)}r, & 0<r\le3/2.
 \end{aligned}
\end{equation}
\end{proposition}
\begin{proof}
Fix a compact interval $J\subset(0,2]$ and choose $b>0$ so that
$[r-b,r+b]\subset(0,4)$ for every $r\in J$. On these intervals $a$ has
a positive lower bound and its first two derivatives are bounded.
For $0<z<b$, the Taylor estimates in the proof of
Lemma~\ref{lem:sharp-small-displacements} therefore give
$N,|\partial_r N|\le C_J z^p$. For $z\ge b$, both increments are
bounded and uniformly separated from zero on $J$. Thus the limit in
\eqref{eq:sharp-paired} exists and equals
\[
 g(r)=\int_0^\infty N(r,z)\frac{\mathrm dz}{z^p}.
\]
Symmetric truncation gives $g(0)=0$ and $g(-r)=-g(r)$.

For $z\ge b$, local absolute continuity of $v$ gives absolute continuity
of $r\mapsto N(r,z)$ on $J$, with
\[
 |\partial_r N(r,z)|\le C_J(1+a(|r-z|))\qquad\text{for a.e. }r\in J.
\]
This bound is integrable on $J\times[b,\infty)$ against
$\mathrm dr\,\mathrm dz/z^p$ since $a\in L^1(0,\infty)$.
By Fubini's theorem, $g$ is absolutely continuous on $J$, with
\[
 g'(r)=\int_0^\infty\partial_r N(r,z)\frac{\mathrm dz}{z^p}
 \qquad\text{for a.e. }r\in J.
\]
The integral over $z<b$ is continuous in $r$ by dominated convergence.
For $z\ge b$, the coefficients $d_\pm^{p-2}$ vary uniformly continuously
with $r$, and the translates of $a(|\cdot|)$ vary continuously in
$L^1(\R)$. Hence $g$ is $C^1$ on $(0,2)$ and the derivative formula
holds there.

Adding the estimates in Lemmas~\ref{lem:sharp-small-displacements}
and~\ref{lem:sharp-large-displacements} gives \eqref{eq:sharp-rhs}
for $0<r<1/16$. On $J=[1/16,2]$, the bounds above are uniform in
$\omega$ by \eqref{eq:sharp-derivatives}, the uniform $L^1$ bound on
$a$ and the fixed choice of $a$ on $[1,\infty)$. This bounds $g$ and
$g'$ on $J$. Since $\Omega(r)\ge\min\{r,1\}$, the derivative bound
in \eqref{eq:sharp-rhs} follows also for $1/16\le r\le3/2$.
\end{proof}

\subsection{The sharpness example and the Dini criterion}\label{subsec:sharp-example}

\begin{proof}[Proof of Theorem~\ref{thm:sharp-modulus}]
We first verify the one-dimensional weak equation. By \eqref{eq:sharp-paired},
$0\le g_\varepsilon(r)\le g(r)$ for $0<r<2$. Oddness gives the same
absolute bound on $(-2,0]$. For $\zeta\in C_c^\infty((-2,2))$ extended
by zero, symmetry gives
\[
 \int_{-2}^2 g_\varepsilon\zeta\,\mathrm dx
 =\frac12\iint_{|x-y|>\varepsilon}
 \frac{J_p(v(x)-v(y))(\zeta(x)-\zeta(y))}{|x-y|^{1+sp}}
 \,\mathrm dx\,\mathrm dy.
\]
The full double integral is absolutely convergent by
$v\in W^{1,p}_{\mathrm{loc}}(\R)\subset W^{s,p}_{\mathrm{loc}}(\R)$
and the boundedness of $v$, as in \eqref{eq:weakform}. Dominated
convergence on both sides proves the weak equation on $(-2,2)$.

Set $u(x)=v(x_1)$ on $\R^n$ and $f(x)=c_{n,p}g(x_1)$ in $B_2$.
Then $u$ is bounded and belongs to $W^{1,p}_{\mathrm{loc}}(\R^n)$,
since $\nabla u=v'(x_1)e_1$ almost everywhere. For $n=1$ take
$c_{1,p}=1$. For $n\ge2$, write $x=(x_1,x')$ and $y=(y_1,y')$.
For $x_1\ne y_1$, integration in the transverse variables gives
\[
 \int_{\R^{n-1}}\frac{\mathrm dz'}{(|x_1-y_1|^2+|z'|^2)^{(n+sp)/2}}
 =\frac{c_{n,p}}{|x_1-y_1|^{1+sp}},\qquad
 c_{n,p}=\int_{\R^{n-1}}\frac{\mathrm dz'}{(1+|z'|^2)^{(n+sp)/2}}.
\]
For $\zeta\in C_c^\infty(B_2)$, the function
$\bar\zeta(x_1)=\int_{\R^{n-1}}\zeta(x_1,x')\,\mathrm dx'$ belongs to
$C_c^\infty((-2,2))$. In the absolutely convergent weak integral for
$u$, fix $x_1\ne y_1$ and integrate the difference
$\zeta(x_1,x')-\zeta(y_1,y')$ against the kernel in $x',y'$.
The result is $c_{n,p}|x_1-y_1|^{-1-sp}
(\bar\zeta(x_1)-\bar\zeta(y_1))$. Fubini's theorem and the
one-dimensional equation tested with $\bar\zeta$ give
\eqref{eq:weakform} with equality in $B_2$. The uniform norm bounds
follow from Lemma~\ref{lem:sharp-profile} and Proposition~\ref{prop:sharp-rhs}.

For the upper bound on $\omega(f;r)$, take $r_*=1/64$ and let $0<r<r_*$ and
$x,y\in B_{3/2}$ with $|x-y|\le r$. If $x_1,y_1$ have opposite signs
or $\min\{|x_1|,|y_1|\}\le2r$, then $|x_1|,|y_1|\le3r$.
By \eqref{eq:sharp-rhs} and concavity,
$|g(x_1)-g(y_1)|\le C\Omega(3r)\le3C\Omega(r)$.
For points of the same sign with $d=\min\{|x_1|,|y_1|\}>2r$, use the
derivative estimate in \eqref{eq:sharp-rhs}. Since $\Omega(t)/t$ is
nonincreasing,
\[
 |g(x_1)-g(y_1)|\le Cr\,\Omega(d)/d\le C\Omega(r).
\]
For the lower bound, the values of $f$ at $(r/4)e_1$ and $-(r/4)e_1$
differ by $2c_{n,p}g(r/4)\ge c\Omega(r)$. Since $g$ is continuous away
from zero, $f(x)=c_{n,p}g(x_1)$ is continuous off $\{x_1=0\}$.
The two points are at distance $r/2$, so sufficiently small
balls about them contain a set of pairs of positive product
measure with distance at most $r$ and difference at least $c\Omega(r)/2$.
Since $\Omega\asymp\omega$, this proves \eqref{eq:sharp-data}. The upper
bound in \eqref{eq:sharp-rhs} gives continuity across $\{x_1=0\}$
when $\omega(0+)=0$.

Finally, $u(re_1)-u(0)=v(r)$ has the order in \eqref{eq:sharp-increment}
by \eqref{eq:sharp-growth}. As $v$ is odd and concave on the positive
half-line, its increment over any interval of length at most $r$ is
at most $2v(r)$. Taking $(re_1,0)$ for the lower bound gives
\[
 v(r)\le\sup\{|u(x)-u(y)|:x,y\in B_1,\ |x-y|\le r\}\le2v(r).
 \qedhere
\]
\end{proof}

\begin{corollary}\label{cor:sharp-dini}
Let $\omega$ satisfy the assumptions of Theorem~\ref{thm:sharp-modulus}.
Among the local weak solutions in Theorem~\ref{thm:logmodulus} satisfying
\[
 \norm{u}_{L^\infty(B_{3/2})}+\Tail(u;0,3/2)
       +\norm{f}_{L^\infty(B_{3/2})}^{1/(p-1)}\le1,
 \qquad \omega(f;t)\le\omega(t)\quad(0<t\le1),
\]
the seminorm $[u]_{\Lip(B_1)}$ is uniformly bounded if and only if
$\int_0^1\omega(t)\,\mathrm dt/t<\infty$.
\end{corollary}
\begin{proof}
Sufficiency follows from Corollary~\ref{cor:dini}. If the integral
diverges, the solution of Theorem~\ref{thm:sharp-modulus} satisfies
\[
 \frac{u(re_1)-u(0)}r\longrightarrow\infty\qquad(r\downarrow0).
\]
For this datum $\omega(f;t)\le C\omega(t)$ on $(0,r_*)$ by
\eqref{eq:sharp-data}. On $[r_*,1]$ the same bound follows from
boundedness of $f$ and $\omega(r_*)>0$, the latter by subadditivity
and $\omega(1)=1$. Multiplying $u$ by a sufficiently small positive
constant and $f$ by its $(p-1)$st power gives both normalizations
while preserving the divergence.
\end{proof}

\Needspace{7\baselineskip}
\begin{remark}[Bounded and logarithmic data]\label{rem:sharp-endpoints}
Taking $\omega(t)=1$ in Theorem~\ref{thm:sharp-modulus} gives the
bounded-data rate $r(1+\log(1/r))^{1/(p-1)}$ in \eqref{eq:logestimate}.
For $\omega(t)=(1-\log t)^{-\mu}$ with $0<\mu\le1$, the constructed
datum is continuous and \eqref{eq:sharp-increment} gives
\[
 u(re_1)-u(0)\asymp
 \begin{cases}
  r(1-\log r)^{(1-\mu)/(p-1)},&0<\mu<1,\\
  r[1+\log(1-\log r)]^{1/(p-1)},&\mu=1.
 \end{cases}
\]
For data with $\omega(f;t)\le C\omega(t)$, the corresponding upper
bounds follow from \eqref{eq:mainmodulus}.
For $\mu>1$ a data modulus bounded by $C(1-\log t)^{-\mu}$ is Dini.
Thus $\mu>1$ is the sharp threshold for a uniform interior Lipschitz
estimate in this logarithmic scale. Corollary~\ref{cor:dini} also covers
slower Dini moduli such as
$(1-\log t)^{-1}[1+\log(1-\log t)]^{-1-\varepsilon}$ for $\varepsilon>0$.
\end{remark}

\section*{Acknowledgments}
\sectionmark{Acknowledgments}
Generative AI tools were used to assist with calculations and checking of parts of the arguments, and the preparation
of the exposition. Responsibility for the mathematical claims and
references rests with the author.

\end{document}